\input ./style/arxiv-general.cfg
\documentclass[aap,MSNbibl,seceqn,dvips]{arximspdf}
\makeatletter
   \@ifpackageloaded{graphicx}{}{\usepackage{graphicx}}
\makeatother
\usepackage{mathbh}

\doi{10.1214/14-AAP1087}
\volume{26}
\issue{1}
\pubyear{2016}
\firstpage{136}
\lastpage{185}
\docsubty{FLA}

\makeatletter

\newcommand{\rrvert}{\vert}
\newcommand{\llvert}{\vert}
\renewcommand{\mid}{|}
\newtheorem{theorem}{Theorem}[section]
\newtheorem{lemma}[theorem]{Lemma}
\newtheorem{proposition}[theorem]{Proposition}
\newproclaim{definition}[theorem]{Definition}
\newproclaim{remark}[theorem]{Remark}
\newproclaim{param}[theorem]{Parameter Optimisation}
\newproclaim{example}[theorem]{Example}
\newcommand{\II}{\mathbb{I}}
\newcommand{\IT}{\mathbb{T}}
\newcommand{\ID}{\mathbb{D}}
\newcommand{\N}{\mathbb{N}}
\newcommand{\Z}{\mathbb{Z}}
\renewcommand{\L}{\mathbb{L}}
\newcommand{\R}{\mathbb{R}}
\renewcommand{\P}{\mathbb{P}}
\newcommand{\cD}{\mathcal{D}}
\newcommand{\cE}{\mathcal{E}}
\newcommand{\cA}{\mathcal{A}}
\newcommand{\cL}{\mathcal{L}}
\newcommand{\cU}{\mathcal{U}}
\newcommand{\cN}{\mathcal{N}}
\newcommand{\cZ}{\mathcal{Z}}
\newcommand{\cY}{\mathcal{Y}}
\newcommand{\cT}{\mathcal{T}}
\newcommand{\cF}{\mathcal{F}}
\newcommand{\sig}{\sigma}
\newcommand{\E}{\mathbb{E}}
\newcommand{\Var}{\operatorname{Var}}
\newcommand{\1}{\mathbh{1}}
\newcommand{\dd}{\mathrm{d}}
\makeatother

\begin{document}
\begin{frontmatter}

\title{Multilevel Monte Carlo for L\'evy-driven SDE\textup{s}: Central limit
theorems for adaptive Euler schemes}
\runtitle{MLMC for adaptive Euler schemes}

\begin{aug}
\author[A]{\fnms{Steffen}~\snm{Dereich}\corref{}\ead[label=e1]{steffen.dereich@wwu.de}}
\and
\author[A]{\fnms{Sangmeng}~\snm{Li}\thanksref{T1}\ead[label=e2]{sangmeng.li@wwu.de}}
\runauthor{S. Dereich and S. Li}
\affiliation{Westf\"alische Wilhelms-Universit\"at M\"unster}
\address[A]{Institute for Mathematical Statistics\\
Westf\"alische Wilhelms-Universit\"at M\"unster\\
Orl\'eans-Ring 10\\
48149 M\"unster\\
Germany\\
\printead{e1}\\
\phantom{E-mail: }\printead*{e2}}
\end{aug}
\thankstext{T1}{Supported by the Deutsche Forschungsgemeinschaft (DFG)
within the Priority
Programme 1324 in the project ``Constructive Quantization and
Multilevel Algorithms for Quadrature of
SDEs.''}

%
\received{\smonth{3} \syear{2014}}
%
\revised{\smonth{11} \syear{2014}}

%
\begin{abstract}
In this article, we consider multilevel Monte Carlo for the numerical
computation
of expectations for stochastic differential equations driven by L\'evy
processes. The underlying numerical schemes are based on jump-adapted
Euler schemes. We prove stable convergence of an idealised scheme.
Further, we deduce limit theorems for certain classes of functionals
depending on the whole trajectory of the process. In particular, we
allow dependence on marginals, integral averages and the supremum of
the process. The idealised scheme is related to two practically
implementable schemes and corresponding\vspace*{1pt} central limit theorems are
given. In all cases, we obtain errors of order $N^{-1/2} (\log
N)^{1/2}$ in the computational time $N$ which is the same order as
obtained in the classical set-up analysed by Giles~[\textit{Oper. Res.} \textbf{56} (2008) 607--617].
Finally, we use the central limit theorems to optimise the parameters
of the multilevel scheme.
\end{abstract}

%
\end{frontmatter}

\section{Introduction}

The numerical computation of expectations $\E[F(X)]$ for solutions
$(X_t)_{t\in[0,T]}$ of stochastic differential equations (SDEs) is a
classical problem in stochastic analysis and numerous numerical schemes
were developed and analysed within the last twenty years; see, for
instance, the textbooks by Kloeden and Platen~\cite{KloedenPlaten} and
Glasserman~\cite{Paul04}. Recently, a new very efficient class of
Monte Carlo algorithms was introduced by Giles~\cite{giles}; see also
Heinrich~\cite{heinrich} for an earlier variant of the computational
concept. Central to these \emph{multilevel} Monte Carlo algorithms is
the use of whole hierarchies of approximations in numerical simulations.
For SDEs, multilevel algorithms often achieve errors of order
$N^{-1/2+o(1)}$ in the computational time $N$ (see \cite
{dereichheidenreich,giles}) despite the infinite-dimensional
nature of the stochastic differential equation. Further, the algorithms
are in many cases optimal in a worst case sense \cite{Creutzigdereich09}.
So far, the main focus of research was concerned with asymptotic error
estimates, whereas central limit theorems have only found minor
attention yet. Beyond the central limit theorem, developed by Ben Alaya
and Kebaier \cite{kebaier} for the Euler scheme for diffusions no
further results are available yet. In general, central limit theorems
illustrate how the choice of parameters affects the efficiency of the
scheme and they are a central tool for tuning the parameters.

%
In this article, we focus on central limit theorems for L\'evy-driven
stochastic differential equations. We prove stable convergence of the
error process of an idealised jump-adapted Euler schemes. Based on this
result, we derive central limit theorems for multilevel schemes for the
approximate computation of expectations of functionals depending on
marginals, integral averages and the supremum of the SDE. We then
introduce implementable jump-adapted Euler schemes that inherit the
properties of the idealised schemes so that the main results prevail.
Finally, we use our new results to optimise over the parameters of the
scheme, and thereby complement the research conducted in~\cite{giles}.
In the Parameter Optimisation~\ref{re1711-1} below, we find that often
it is preferable to increase the number of Euler steps from level to
level by a factor of $6$.
For ease of presentation, we restrict attention to the one-dimensional
setting although a generalisation to finite-dimensional stochastic
differential equations is canonical.

In the following, $(\Omega,\cF,\P)$ denotes a probability space that
is sufficiently rich to ensure existence of all random variables used
in the exposition. We let $Y=(Y_t)_{t\in[0,T]}$ be a square integrable
L\'evy-process and note that there exist $b\in\R$ (drift), $\sigma
^2\in[0, \infty)$ (diffusion coefficient) and a measure $\nu$ on $\R
\setminus\{0\}$ with $\int x^2 \nu(\dd x)<\infty$ (L\'evy measure)
such that
\[
\E\bigl[e^{i z Y_t}\bigr]= \exp \biggl\{ t \biggl( ibz-\frac{1}{2}
\sigma^2 z + \int\bigl(e^{izx}-1-izx\bigr) \nu(\dd x) \biggr)
\biggr\} %
\]
for $t\in[0,T]$ and $z\in\R$. We call the unique triplet $(b,\sig
^2,\nu)$ the L\'evy triplet of $Y$. We refer the reader to the
textbooks by Applebaum~\cite{appl}, Bertoin~\cite{bert} and
Sato~\cite{sato} for a concise treatment of L\'evy processes.
The process $X=(X_t)_{t\in[0,T]}$ denotes the solution to the
stochastic integral equation
%
\begin{equation}
\label{sde} X_{t}=x_{0}+\int_{0}^{t}a(X_{s-})\,\dd Y_{s}, \qquad t\in[0,T],
\end{equation}
where $a\dvtx \R\to\R$ is a continuously differentiable Lipschitz
function and $x_0\in\R$. Both processes $Y$ and $X$ attain values in
the space of c\`adl\`ag functions on $[0,T]$ which we will denote by
$\ID(\R)$ and endow with the Skorokhod topology.
We will analyse multilevel algorithms for the computation of expectations
$\E[F(X)]$,
where $F\dvtx \ID(\R)\to\R$ is a measurable functional such that $F(x)$
depends on the marginals, integrals and/or supremum of the path $x\in
\ID(\R)$.
Before we state the results, we introduce the underlying numerical schemes.

\subsection{Jump-adapted Euler scheme}\label{sec1803-02}

In the context of L\'evy-driven stochastic differential equations,
there are various Euler-type schemes analysed in the literature. We
consider jump-adapted Euler schemes. For finite L\'evy measures, these
were introduced by Platen~\cite{platen} and analysed by various
authors; see, for example, \mbox{\cite{Maghsoodi,bruti}}. For infinite L\'evy
measures, an error analysis is conducted in \cite{dereichheidenreich}
and \cite{dereich2} for two multilevel Monte Carlo schemes. Further,
weak approximation is analysed in \cite{kohatsutankov} and \cite
{mordecki}. 
In general, the simulation of increments of the L\'evy-process is
delicate. One can use truncated shot noise representations as in~\cite
{rosinski01}. These perform well for Blumenthal--Getoor indices smaller
than one, but are less efficient when the BG-index gets larger than
one~\cite{dereichheidenreich}, even when combined with a Gaussian
compensation in the spirit of~\cite{sorenjan}; see~\cite{Dereich}. A
faster simulation technique is to do an inversion of the characteristic
function of the L\'evy process and to establish direct simulation
routines in a precomputation. Certainly, this approach is more
involved and its realisation imposes severe restrictions on the
dimension of the L\'evy process; see~\cite{DerLi14a}.


In this article, we analyse one prototype of adaptive approximations
that is intimately related to implementable adaptive schemes and we
thus believe that our results have a universal appeal. The
approximations depend on two positive parameters:
\begin{itemize}
\item$h$, the threshold for the size of the jumps being considered
large and causing immediate updates, and
\item$\varepsilon$ with $T\in\varepsilon\N$, the length of the
regular update intervals.
\end{itemize}
For the definition of the approximations, we use the simple Poisson
point process~$\Pi$ on the Borel sets of $(0,T]\times(\R\setminus\{
0\})$ associated to $Y$, that is,
\[
\Pi= \sum_{s\in(0,T]\dvtx  \Delta Y_s \neq0} \delta_{(s,\Delta Y_s)}, %
\]
where we use the notation $\Delta x_t=x_t-x_{t-}$ for $x\in\ID(\R)$
and $t\in(0,T]$.
It has intensity $\ell_{(0,T]}\otimes\nu$, where $\ell_{(0,T]}$
denotes Lebesgue measure on $(0,T]$. Further, let $\overline\Pi$ be the
compensated variant of $\Pi$ that is the random signed measure on
$(0,T]\times(\R\setminus\{0\})$ given by
\[
\overline\Pi=\Pi- \ell_{(0,T]}\otimes\nu. %
\]
The process $(Y_t)_{t\in[0,T]}$ admits the representation
%
\begin{equation}
\label{firstdecom} Y_{t}=bt+\sigma W_{t}+ \lim
_{\delta\downarrow0} \int_{(0,t]\times
B(0,\delta)^c } x\,\dd\overline{\Pi}(s,x),
\end{equation}
where $(W_t)_{t\in[0,T]}$ is an appropriate Brownian motion that is
independent of $\Pi$ and the limit is to be understood uniformly in
$\mathbb{L}^{2}$. 
We enumerate the random set
\[
\bigl(\varepsilon\Z\cap[0,T] \bigr) \cup\bigl\{t\in(0,T]\dvtx  \llvert \Delta
Y_t\rrvert \geq h\bigr\}=\{T_0,T_1,\ldots\}
\]
in increasing order and define the approximation $ X^{h,\varepsilon}
=( X^{h,\varepsilon}_t)_{t\in[0,T]}$ by $X^{h,\varepsilon}_0=x_0$
and, for $n=1,2,\ldots$ and $t\in(T_{n-1},T_{n}]$
%
\begin{equation}
\label{update} X^{h,\varepsilon}_{t}=X^{h,\varepsilon}_{T_{n-1}}+a
\bigl( X^{h,\varepsilon
}_{T_{n-1}}\bigr) (Y_t-Y_{T_{n-1}}).
\end{equation}

\subsection{Multilevel Monte Carlo}\label{secMLintro}

In general, multilevel schemes make use of whole hierarchies of
approximate solutions and we choose decreasing sequences $(\varepsilon
_k)_{k\in\N}$ and $(h_k)_{k\in\N}$ with:
\begin{longlist}[(ML2)]
\item[(ML1)] $\varepsilon_k=M^{-k} T$, where $M\in\{2,3,\ldots\}$ is fixed,

\item[(ML2)] $\lim_{k\to\infty} \nu(B(0,h_k)^c) \varepsilon_k=
\theta$ for a $\theta\in[0,\infty)$ and
$\lim_{k\to\infty} h_k/\sqrt{\varepsilon_k}=0$.
\end{longlist}

We remark that whenever $\theta$ in \textup{(ML2)} is strictly positive, then
one automatically has that $h_k=o(\sqrt{\varepsilon_k})$; see
Lemma~\ref{driftdisappear}.

For every $k\in\N$, we denote by $X^k:=X^{h_k,\varepsilon_k}$
the corresponding adaptive Euler approximation with update rule~(\ref
{update}). Once\vspace*{1pt} this hierarchy of approximations has been fixed, a multilevel
scheme $\widehat S$ is parameterised by a $\N$-valued vector
$(n_1,\ldots,n_L)$ of arbitrary finite length $L$: for a measurable
function $F\dvtx \ID(\R)\to\R$ we approximate $\E[F(X)]$ by
\[
\E\bigl[F\bigl( X^1\bigr)\bigr]+ \E\bigl[F\bigl(X^2
\bigr)-F\bigl(X^{1}\bigr)\bigr]+ \cdots+ \E\bigl[F\bigl(
X^L\bigr)-F\bigl( X^{L-1}\bigr)\bigr] %
\]
and denote by $\widehat S(F)$ the random output that is obtained when
estimating the individual expectations $\E[F(X^1)],  \E[F(
X^2)-F(X^{1})],  \ldots,  \E[F( X^L)-F(X^{L-1})]$ independently by
classical Monte Carlo with $n_1,\ldots,n_L$ iterations and summing up
the individual estimates. More explicitly, a multilevel scheme
$\widehat S$ associates to each measurable $F$ a random variable
%
\begin{equation}
\label{multiestimator} \widehat{S}(F)=\frac{1}{n_{1}}\sum_{i=1}^{n_{1}}F
\bigl( X^{1,i}\bigr)+\sum_{k=2}^{L}
\frac{1}{n_{k}}\sum_{i=1}^{n_{k}} \bigl(F
\bigl( X^{k,i,f}\bigr)-F\bigl( X^{k-1,i,c}\bigr) \bigr),
\end{equation}
where the pairs of random variables $( X^{k,i,f}, X^{k-1,i,c})$,
respectively, the random variables $ X^{1,i}$, appearing in the sums
are all independent with identical distribution as $(X^k,X^{k-1})$,
respectively,~$X^1$. Note that the upper indices $f$ and $c$ refer to
fine and coarse and that the entries of each pair are \emph{not} independent.

\subsection{Implementable schemes}

We give two implementable schemes. The first one relies on
precomputation for direct simulation of L\'evy increments. The second
one ignores jumps of size smaller than a threshhold which leads to
schemes of optimal order only in the case where---roughly speaking---the Blumenthal--Getoor index is smaller than one.

\subsubsection*{Schemes with direct simulation of small jumps}

For $h>0$, we let
$Y^h=(Y^h_t)_{t\in[0,T]}$ denote the L\'evy process given by
%
\begin{equation}
\label{defYh} Y^h_t= bt +\sigma W_t + \int
_{(0,t]\times B(0,h)^c } x \,\dd\overline{\Pi}(s,x).
\end{equation}
Using\vspace*{1pt} the shot noise representation (see \cite{sorenjan}), we can
simulate $Y^h$ on arbitrary (random) time sets. The remainder
$M^h=(M^h_t)_{t\in[0,T]}$, that is,
\[
M_t^h=\lim_{\delta\downarrow0} \int
_{(0,t]\times(B(0,h)\setminus
B(0,\delta))} x \,\dd\overline\Pi(s,x)= Y-Y^h, %
\]
can be simulated on a fixed time grid $\varepsilon' \Z\cap[0,T]$
with $\varepsilon'\in\varepsilon\N$ denoting an additional
parameter of the scheme.
A corresponding approximation is given by $\widehat{X}^{h,\varepsilon,\varepsilon'}=(\widehat{X}^{h,\varepsilon,\varepsilon'}_t)_{t\in[0,T]}$
via $\widehat{X}^{h,\varepsilon,\varepsilon'}_0=x_0$ and, for $n=1,2,\ldots
$ and $t\in(T_{n-1},T_n]$,
%
\begin{eqnarray}\label{continuousapp}
\widehat{X}^{h,\varepsilon,\varepsilon'}_t&=&\widehat{X}^{h,\varepsilon,\varepsilon'}_{T_{n-1}}+a
\bigl(\widehat{X}^{h,\varepsilon,\varepsilon
'}_{T_{n-1}}\bigr) \bigl(Y^h_{t}-Y^h_{T_{n-1}}
\bigr)
\nonumber\\[-8pt]\\[-8pt]\nonumber
&&{} + \1_{\varepsilon'\Z}(t) a\bigl(\widehat{X}^{h,\varepsilon,\varepsilon'}_{t-\varepsilon'}\bigr)
\bigl(M^h_{t}-M^h_{t-\varepsilon'}\bigr).
\end{eqnarray}
We call $\widehat{X}^{h,\varepsilon,\varepsilon'}$ the \emph{continuous
approximation} with parameters $h,\varepsilon,\varepsilon'$. Further,
we define the \emph{piecewise constant approximation} $\overline X^{h,\varepsilon,\varepsilon'}=(\overline X^{h,\varepsilon,\varepsilon
'}_t)_{t\in[0,T]}$ via demanding that for $n=1,2,\ldots$ and $t\in
[T_{n-1},T_n)$,
%
\begin{equation}
\label{constantapp} \overline X^{h,\varepsilon,\varepsilon'}_{t}=\widehat{X}^{h,\varepsilon,\varepsilon'}_{T_{n-1}}
\end{equation}
and $\overline X^{h,\varepsilon,\varepsilon'}_{T}=\widehat{X}^{h,\varepsilon,\varepsilon'}_{T}$.

In corresponding multilevel schemes, we choose $(\varepsilon_k)_{k\in
\N}$ and $(h_k)_{k\in\N}$ as before. Further, we choose
monotonically decreasing parameters $(\varepsilon'_k)_{k\in\N}$ with
$\varepsilon_k'\in\varepsilon_k\N$ and:
\begin{longlist}[(ML3a)]
\item[(ML3a)] $\varepsilon_k'\int_{B(0,h_k)}x^2 \nu(\dd x)\log^2
(1+1/\varepsilon_k')=o(\varepsilon_k)$,
\item[(ML3b)] $h_k^2\log^2 (1+1/\varepsilon_k')=o(\varepsilon_k)$.
\end{longlist}

\begin{remark}
If
%
\begin{equation}
\label{eq2602-3} \int x^2 \log^2 \biggl(1+\frac{1}x
\biggr) \nu(\dd x)<\infty,
\end{equation}
there exist appropriate parameters $(h_k,\varepsilon_k,\varepsilon
'_k)_{k\in\N}$ satisfying \textup{(ML1)}, \textup{(ML2)}, \textup{(ML3a)} and \textup{(ML3b)}.
More precisely, in the case where $\nu$ is infinite, appro\-priate
parameters are obtained by choosing $\varepsilon'_k=\varepsilon_k$
and $(h_k)$ with\break $\lim_{k\to\infty} \varepsilon_k  \nu
(B(0,h_k)^c)\hspace*{-0.6pt}=\theta>0$; see Lemma~\ref{driftdisappear}.
\end{remark}

In analogy to before, we denote by $(\widehat{X}^k\dvtx k\in\N)$ and $(\overline X^k\dvtx k\in\N)$ the corresponding approximate continuous and piecewise
constant solutions. We state a result of~\cite{Li15} which implies
that in most cases the central limit theorems to be provided later are
also valid for the continuous approximations.

\begin{lemma}\label{le2602-1}If assumptions \textup{(ML1)}, \textup{(ML3a)} and \textup{(ML3b)}
are satisfied, then
%
\[
\lim_{k\to\infty} \varepsilon_k^{-1} \E \Bigl[
\sup_{t\in[0,T]} \bigl\llvert X^k_t-\widehat{X}^k_t\bigr\rrvert ^2 \Bigr]=0. %
\]
\end{lemma}

Practical issues of numerical schemes with direct simulation of
increments are discussed in~\cite{DerLi14a}. 

\subsubsection*{Truncated shot noise scheme}

The truncated shot noise scheme is parameterised by two positive
parameters $h,\varepsilon$ as above.
The \emph{continuous approximations} $\widehat{X}^{h,\varepsilon}=(\widehat{X}^{h,\varepsilon}_t)_{t\in[0,T]}$ are defined via $\widehat{X}^{h,\varepsilon}_0=x_0$ and, for $n=1,2,\ldots$ and $t\in(T_{n-1},T_n]$,
%
\begin{equation}
\label{continuousapp2} \widehat{X}^{h,\varepsilon}_t=\widehat{X}^{h,\varepsilon}_{T_{n-1}}+a
\bigl(\widehat{X}^{h,\varepsilon}_{T_{n-1}}\bigr) \bigl(Y^h_{t}-Y^h_{T_{n-1}}
\bigr)
\end{equation}
and the \emph{piecewise constant approximations} $\overline X^{h,\varepsilon}=(\overline X^{h,\varepsilon}_t)_{t\in[0,T]}$ are
defined as before by
demanding that, for $n=1,2,\ldots$ and $t\in[T_{n-1},T_n)$,
%
\begin{equation}
\label{constantapp2} \overline X^{h,\varepsilon}_{t}=\widehat{X}^{h,\varepsilon}_{T_{n-1}}
\end{equation}
and $\overline X^{h,\varepsilon}_{T}=\widehat{X}^{h,\varepsilon}_{T}$. Again we
will use decreasing sequences $(\varepsilon_k)$ and $(h_k)$ as before
to specify sequences of approximations $(\widehat{X}^k)$ and $(\overline X^k)$.
In the context of truncated shot noise schemes, we will impose as
additional assumption:
\begin{longlist}[(ML4)]
\item[(ML4)] $\int_{B(0,h_k)}x^2 \nu(\dd x)=o(\varepsilon_k)$.
\end{longlist}

\begin{remark}
If $\int\llvert  x\rrvert   \nu(\dd x)<\infty$, then \textup{(ML1)}, \textup{(ML2)} and \textup{(ML4)} are
satisfied for appropriate parameters.
\end{remark}

The following result is a minor modification of \cite{dereichheidenreich}, Proposition~1; see also~\cite{Li15}.

\begin{lemma}\label{le2602-2} If assumptions \textup{(ML1)} and \textup{(ML4)} are
satisfied, then
\[
\lim_{k\to\infty} \varepsilon_k^{-1} \E \Bigl[
\sup_{t\in[0,T]} \bigl\llvert X^k_t-\widehat{X}^k_t\bigr\rrvert ^2 \Bigr]=0. %
\]
\end{lemma}

\subsection{Main results}

In the following, we will always assume that $Y\hspace*{-0.2pt}=(Y_t)_{t\in[0,T]}$ is
a square integrable L\'evy process with L\'evy triplet $(b,\sig^2,\nu
)$ satisfying $\sig^2>0$ and that $X=(X_t)_{t\in[0,T]}$ solves the SDE
\[
\dd X_t=a(X_{t-}) \,\dd Y_t %
\]
with $X_0=x_0$, where $a\dvtx \R\to\R$ is a continuously differentiable
Lipschitz function. Further, for each $k\in\N$, $X^k$ denotes the
jump-adapted Euler scheme with updates at all times in
\[
\bigl(\varepsilon_k \N\cap[0,T]\bigr)\cup\bigl\{t\in(0,T]\dvtx  \llvert
\Delta Y_t\rrvert \geq h_k\bigr\}; %
\]
see~(\ref{update}). The decreasing sequences of parameters
$(\varepsilon_k)$ and $(h_k)$ are assumed to satisfy (ML1) and \textup{(ML2)}
from Section~\ref{secMLintro}.%

\subsubsection*{Convergence of the error process}
We consider the normalised sequence of \emph{error processes}
associated to the multilevel scheme that is the sequence $(\varepsilon
_k^{-1/2}(X^{k+1}-X^k)\dvtx k\in\N)$. Let us introduce the process
appearing as a limit.
We equip the points of the associated point process $\Pi$ with
independent marks and denote for a point $(s,x)\in\Pi$:
\begin{itemize}
%
\item by $\xi_{s}$, a standard normal random variable,
\item by $\mathcal U_{s}$, an\vspace*{1pt} independent uniform random variable on
$[0,1]$, and
\item by $\mathcal E^{\theta}_{s}$ and $\mathcal E^{(M-1)\theta}_{s}$
independent $\operatorname{Exp}(\theta)$ and $\operatorname{Exp}((M-1)\theta)$-distributed random
variables, respectively.
\end{itemize}
Further, we denote by $B=(B_t)_{t\in[0,T]}$ an independent standard
Brownian motion.

The \emph{idealised error process} $U=(U_{t})_{t\in[0,T]}$ is defined
as the solution of the integral equation
%
\begin{eqnarray}
\label{limitproce} U_{t}&=&\int_{0}^{t}a'(X_{s-})U_{s-}
\,\dd Y_{s}+\sigma^2 \Upsilon \int_{0}^{t}
\bigl(aa'\bigr) (X_{s-}) \,\dd B_{s}
\nonumber\\[-8pt]\\[-8pt]\nonumber
&&{}+\sum_{s\in(0,t]\dvtx  \Delta Y_s\neq0}\sigma_s \mathcal
\xi_{s} \bigl(aa'\bigr) (X_{s-}) \Delta
Y_{s},
\end{eqnarray}
where $\Upsilon^2=\frac{e^{-\theta}-1+\theta}{\theta^2} (1-\frac{1}M)$, if $\theta>0$, and $\Upsilon^2=\frac{1}2(1-\frac{1}{M})$, if
$\theta=0$, and the positive marks $(\sigma_s)$ are defined by
\begin{eqnarray*}
\sigma_s^2&=&\sigma^2\sum
_{1\leq m\leq M}\1_{\{(m-1)/{M}\leq\cU
_{s}<m/M\}} \biggl[\min\bigl(
\cE_{s}^{\theta},\cU_{s}\bigr)
\\
&&\hspace*{138pt}{} -\min\biggl(\cE
_{s}^\theta,\cE_s^{(M-1)\theta},
\cU_{s}-\frac{m-1}{M}\biggr) \biggr].
\end{eqnarray*}
Note that the above infinite sum has to be understood as an appropriate
martingale limit. More explicitly, denoting by $L=(L_t)_{t\in[0,T]}$
the L\'evy process
\[
L_t=\sigma^2 \Upsilon B_{t}+ \lim
_{\delta\downarrow0} \sum_{s\in(0,t]\dvtx  \llvert  \Delta Y_s\rrvert  \geq\delta}
\sigma_s \mathcal\xi_{s} \Delta Y_{s}
\]
we can rewrite~(\ref{limitproce}) as
\[
U_t=\int_{0}^{t}a'(X_{s-})U_{s-}
\,\dd Y_{s}+\int_0^t
\bigl(aa'\bigr) (X_{s-}) \,\dd L_{s}.
\]
Strong uniqueness and existence of the solution follow from Jacod and
Memin \cite{JacMem81}, Theorem~4.5.

\begin{theorem}\label{twosuccessive}
Under the above assumptions, we have weak convergence
%
\begin{equation}
\bigl(Y,\varepsilon_{n}^{-1/2}\bigl( X^{n+1}-
X^{n}\bigr)\bigr) \Rightarrow(Y, U)\qquad\mbox{in }\ID\bigl(
\R^2\bigr).
\end{equation}
\end{theorem}

\subsubsection*{Central limit theorem for linear functionals}

We consider functionals $F\dvtx \ID(\R)\to\R$ of the form
\[
F(x)=f(Ax) %
\]
with $f\dvtx \R^d\to\R$ and $A\dvtx \ID(\R)\to\R^d$ being linear and
measurable. We set
\[
D_f:=\bigl\{z\in\R^d\dvtx  f\mbox{ is differentiable in }z
\bigr\}. %
\]

\begin{theorem}\label{theo2402-1}
Suppose that $f$ is Lipschitz continuous and that $A$ is Lipschitz
continuous with respect to supremum norm and continuous with respect to
the Skorokhod topology in $\P_U$-almost every path. Further suppose that
$AX \in D_f$, almost surely, and that $\alpha\geq\frac{1}2$ is such
that the limit
\[
\lim_{n\rightarrow\infty}\varepsilon^{-\alpha}_n \E
\bigl[F\bigl(X^{n}\bigr)-F(X) \bigr]=:\kappa %
\]
exists.
We denote for $\delta\in(0,1)$ by $\widehat S_\delta$ the multilevel
Monte Carlo scheme with parameters $(n_1^{(\delta)},n_2^{(\delta
)},\ldots,n_{L(\delta)}^{(\delta)})$, where
%
\begin{equation}
\label{eq1103-2} L(\delta)= \biggl\lceil\frac{\log\delta^{-1}}{\alpha\log M} \biggr\rceil\quad
\mbox{and}\quad n_{k}(\delta)= \bigl\lceil\delta ^{-2} L(
\delta) \varepsilon_{k-1} \bigr\rceil,
\end{equation}
for $k=1,2,\ldots, L(\delta)$.
Then we have,
\[
\delta^{-1} \bigl(\widehat S_{\delta}(F)-\E\bigl[F(X)\bigr]
\bigr)\Rightarrow\mathcal {N}\bigl(\kappa,\rho^{2}\bigr)\qquad\mbox{as
}\delta\rightarrow0,
\]
where $\mathcal{N}(\kappa,\rho^{2})$ is the normal distribution with
mean $\kappa$ and variance
\[
\rho^{2}=\Var \bigl( \nabla f(AX) \cdot AU \bigr).
\]
\end{theorem}

\begin{example}
(a)
For any finite signed measure $\mu$, the integral $Ax=\int_0^Tx_s
\,\dd\mu(s)$ satisfies the assumptions of the theorem. Indeed, for
every path $x\in\ID(\R)$ with
%
\begin{equation}
\label{eq0303-2} \mu\bigl(\bigl\{s\in[0,T]\dvtx  \Delta x_s\neq0\bigr\}
\bigr)=0
\end{equation}
one has for $x^n\to x$ in the Skorokhod space that
\[
Ax^n =\int_0^Tx^n_s
\,\dd\mu(s)\to\int_0^Tx_s \,\dd
\mu(s)=Ax %
\]
by dominated convergence and~(\ref{eq0303-2}) is true for $\P
_U$-almost all paths since $\mu$ has at most countably many atoms.
Hence, the linear maps $A x=x_t$ and $Ax=\int_0^T x_s \,\dd s$ are
allowed choices in Theorem~\ref{theo2402-1} since $U$ is almost surely
continuous in $t$.

(b)~All combinations of admissible linear maps $A_1,\ldots,A_m$
satisfy again the assumptions of the theorem.
\end{example}

In view of implementable schemes, we state a further version of the theorem.

\begin{theorem}\label{thm2602-1}
Suppose that either $(\widehat{X}^k\dvtx k\in\N)$ and $(\overline X^k\dvtx k\in\N)$
denote the continuous and piecewise constant approximations of the
scheme with direct simulation and that \textup{(ML1)}, \textup{(ML2)} and \textup{(ML3)} are
fulfilled or that they are the approximations of the truncated shot
noise scheme and that \textup{(ML1)}, \textup{(ML2)} and \textup{(ML4)} are fulfilled. Then\vspace*{1pt}
Theorem~\ref{theo2402-1} remains true when replacing the family
$(X^k\dvtx k\in\N)$ by $(\widehat{X}^k\dvtx k\in\N)$. Further, if $A$ is given by
\[
Ax= \biggl(x_T,\int_0^T
x_s \,\dd s \biggr), %
\]
the statement of the central limit theorem remains true, when replacing
the family $(X^k\dvtx k\in\N)$ by $(\overline X^k\dvtx k\in\N)$.
\end{theorem}

\subsubsection*{Central limit theorem for supremum-dependent functionals}

In this section, we consider functionals $F\dvtx \ID(\R)\to\R$ of the form
\[
F(x)=f \Bigl(\sup_{t\in[0,T]} x_t \Bigr) %
\]
with $f\dvtx \R\to\R$ measurable.

\begin{theorem}\label{theo2402-2}Suppose that $f\dvtx \R\to\R$ is
Lipschitz continuous and that the coefficient~$a$ does not attain zero.
Further, suppose that
$\sup_{t\in[0,T]}X_t \in D_f$, almost surely, and that $\alpha\geq
\frac{1}2$ is such that the limit
\[
\lim_{n\rightarrow\infty}\varepsilon^{-\alpha}_n \E
\bigl[F\bigl(X^{n}\bigr)-F(X) \bigr]=:\kappa %
\]
exists.
We denote for $\delta\in(0,1)$ by $\widehat S_\delta$ the multilevel
Monte Carlo scheme with parameters $(n_1^{(\delta)},n_2^{(\delta
)},\ldots,n_{L(\delta)}^{(\delta)})$, where
\[
L(\delta)= \biggl\lceil\frac{\log\delta^{-1}}{\alpha\log M} \biggr\rceil\quad\mbox{and}\quad
n_{k}(\delta)= \bigl\lceil\delta ^{-2} L(\delta)
\varepsilon_{k-1} \bigr\rceil,
\]
for $k=1,2,\ldots, L(\delta)$.
Then we have
\[
\delta^{-1} \bigl(\widehat S_{\delta}(F)-\E\bigl[F(X)\bigr]
\bigr)\Rightarrow\mathcal {N}\bigl(\kappa,\rho^{2}\bigr)\qquad\mbox{as
}\delta\rightarrow0,
\]
where $\mathcal{N}(\kappa,\rho^{2})$ is the normal distribution with
mean $\kappa$ and variance
\[
\rho^{2}=\Var \Bigl( f' \Bigl(\sup_{t\in[0,T]}
X_t \Bigr) U_S \Bigr),
\]
and $S$ denotes the random time at which $X$ attains its supremum.
\end{theorem}

\begin{theorem}\label{theo2003-1}Theorem~\ref{theo2402-2} remains
true for the continuous approximations for the scheme with direct
simulation of increments or the truncated shot noise scheme under the
same assumptions as imposed in Theorem~\ref{thm2602-1}.
\end{theorem}

\subsubsection*{Optimal parameters}

We use the central limit theorems to adjust the parameters of the
multilevel scheme. Here, we use the following result.

\begin{theorem}\label{theo0303-1}
Let $F$ be as in Theorems~\ref{theo2402-1} or~\ref{theo2402-2} and
assume that the assumptions of the respective theorem are fulfilled.
Further assume in the first case that $A$ is of integral type meaning
that there exist finite signed measures $\mu_1,\ldots,\mu_d$ on
$[0,T]$ such that $A=(A_1,\ldots,A_d)$ with
\[
A_jx=\int_0^Tx_s
\,\dd\mu_j(s)\qquad\mbox{for }x\in\ID(\R)\mbox{ and } j=1,
\ldots,d %
\]
and generally suppose that $a'(X_{s-})\Delta Y_s\neq-1$ for all $s\in
[0,T]$, almost surely.
Then there exists a constant $\kappa$ depending on $F$ and the
underlying SDE, but not on $M$ and $\theta$ such that the variance
$\rho^2$ appearing as variance is of the form
\[
\rho= \kappa \Upsilon, %
\]
where as before $\Upsilon^2=\frac{e^{-\theta}-1+\theta}{\theta^2}
(1-\frac{1}M)$, if $\theta>0$, and $\Upsilon^2=\frac{1}2 (1-\frac{1}M)$,
if $\theta=0$.
\end{theorem}

\begin{remark}
The assumption that $a'(X_{s-})\Delta Y_s\neq-1$ for all $s\in[0,T]$,
almost surely, is automatically fulfilled if $\nu$ has no atoms.
For every $s\in(0,T]$ with $a'(X_{s-})\Delta Y_s=-1$, the error
process jumps to zero causing technical difficulties in our proofs.
In general, the result remains true without this assumption, but for
simplicity we only provide a proof under this technical assumption.
\end{remark}

\begin{figure}

\includegraphics{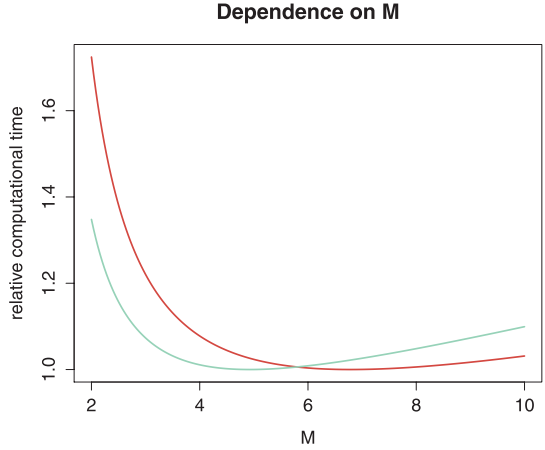}

\caption{Impact of $M$ on the computational cost for $\beta=0$
(green) and $\beta=1$ (red).}\label{fig1}
\end{figure}

\begin{param}\label{re1711-1}
We use Theorem~\ref{theo0303-1} to optimise the parameters. We assume
that $\theta$ of \textup{(ML2)} and the bias $\kappa$ are zero.
Multilevel schemes are based on iterated sampling of
$F(X^{k})-F(X^{k-1})$, where $(X^{k-1},X^{k})$ are coupled approximate
solutions.
Typically, one simulation causes cost (has runtime) of order
\[
C_k=\bigl(1+o(1)\bigr) \kappa_\mathrm{cost}
\varepsilon_{k-1}^{-1} (M+ \beta), %
\]
where $\kappa_\mathrm{cost}$ is a constant that does not depend on
$M$, and $\beta\in\R$ is an appropriate constant typically with values
between zero and one: one coupled path simulation needs:
\begin{itemize}
\item to simulate $\varepsilon_{k-1}^{-1}TM$ increments of the L\'evy process,
\item to do $\varepsilon_{k-1}^{-1}TM$ Euler steps to gain the fine
approximation,
\item to concatenate $\varepsilon_{k-1}^{-1}T(M-1)$ L\'evy increments, and
\item to do $\varepsilon_{k-1}^{-1}T$ Euler steps to gain the coarse
approximation.
\end{itemize}
If every operation causes the same computational cost, one ends up with
$\beta=0$. If the concatenation procedure is significantly less
expensive, the parameter $\beta$ rises.
Using that
\[
\delta^{-1} \bigl(\widehat S_{\delta}(F)-\E\bigl[F(X)\bigr]
\bigr)\Rightarrow\mathcal {N}\bigl(0,\kappa_\mathrm{err}
^{2}(1-1/M)\bigr)\qquad\mbox{as }\delta \downarrow0,
\]
we conclude that for $\bar\delta:=\bar\delta(\delta):=\delta
/(\kappa_\mathrm{err}\sqrt{1-1/M})$ one has
\[
\delta^{-1} \bigl(\widehat S_{\bar\delta}(F)-\E\bigl[F(X)\bigr]
\bigr)\Rightarrow \mathcal{N}(0,1)\qquad\mbox{as }\delta\downarrow0.
\]
Hence, the asymptotics of $\widehat S_{\bar\delta}(F)$ do not depend
on the choice of $M$ and we can compare the efficiency of different
choices of $M$ by looking at the cost of a simulation of $\widehat S_{\bar\delta}(F)$. It is of order
\begin{eqnarray*}
&& \bigl(1+o(1)\bigr) \kappa_{\mathrm{cost}} L(\bar\delta)^2 (M+\beta)
\bar\delta^{-2}
\\
&&\qquad = \bigl(1+o(1)\bigr)\frac{\kappa_{\mathrm{cost}}\kappa_\mathrm
{err}^2}{\alpha^2}
\frac{(M-1)(M+\beta)}{M (\log M)^2} \delta ^{-2} \bigl(\log\delta^{-1}
\bigr)^2. %
\end{eqnarray*}
A plot illustrating the dependence on the choice of $M$ is provided in
Figure~\ref{fig1}. There we plot the function $M\mapsto\frac
{(M-1)(M+\beta)}{M (\log M)^2}$ for $\beta$ being\vspace*{1pt} $0$ or $1$. The
plot indicates that in both cases $6$ is a good choice for $M$. In
particular, it is not necessary to know $\beta$ explicitly in order to
find a ``good'' $M$. For numerical tests concerning appropriate choices
of $\beta$, we refer the reader to~\cite{DerLi14a}.
\end{param}

The article is outlined as follows. In Section~\ref{section1103-1}, we
analyse the error process and prove Theorem~\ref{twosuccessive}. In
Section~\ref{derivedquantities}, we prepare the proofs of the central
limit theorems for integral averages for the piecewise constant
approximations and for supremum dependent functionals. In Section~\ref
{proofCLT}, we provide the proofs of all remaining theorems, in
particular, of all central limit theorems. The article ends with an
\hyperref[append]{Appendix} where we summarise known and auxiliary results. In
particular, we provide a brief introduction to stable convergence and
perturbation estimates mainly developed in articles by Jacod and Protter.


\section{The error process (Theorem~\texorpdfstring{\protect\ref{twosuccessive}}{1.5})}\label{section1103-1}

In this section, we prove Theorem~\ref{twosuccessive}.
We assume that properties~\textup{(ML1)} and~\textup{(ML2)} are fulfilled.
At first, we introduce the necessary notation and outline our strategy
of proof. All intermediate results will be stated as propositions and
their proofs are deferred to later subsections.
We denote for $n\in\N$ and $t\in[0,T]$
\[
\iota_n(t)= \sup [0,t]\cap\II_n, %
\]
where $\II_n=\{s\in(0,T]\dvtx  \Delta Y^{h_n}_s\neq0\}\cup(\varepsilon
_n\Z\cap[0,T])$ is the random set of update times and recall that
$X^n$ solves
%
\begin{equation}
\label{widehatx} \dd X^{n}_t=a\bigl( X^{n}_{\iota_n(t-)}
\bigr) \,\dd Y_t
\end{equation}
with $X^n_0=x_0$.
We analyse the (normalised) \emph{error process} of two consecutive
\mbox{$X^n$-}levels that is the process $U^{n,n+1}=(U^{n,n+1}_t)_{t\in[0,T]}$
given by
\[
U^{n,n+1}_t= \varepsilon_n^{-1/2} \bigl(
X^{n+1}- X^{n}\bigr). %
\]

\noindent
The error process satisfies the SDE
\begin{eqnarray*}
\dd U^{n,n+1}_t&=&\varepsilon_n^{-1/2}
\bigl(a\bigl(X^{n+1}_{t-}\bigr)-a\bigl(X^{n}_{t-}
\bigr) \bigr)\,\dd Y_t+\varepsilon_n^{-1/2} \bigl(a
\bigl( X^{n}_{t-}\bigr)-a\bigl(X^{n}_{\iota_n(t-)}
\bigr) \bigr)\,\dd Y_t
\\
&&{} -\varepsilon_n^{-1/2} \bigl(a\bigl(X^{n+1}_{t-}
\bigr)-a\bigl( X^{n+1}_{\iota
_{n+1}(t-)}\bigr) \bigr)\,\dd Y_t.
\end{eqnarray*}
In order to rewrite the SDE, we introduce some more notation. We let
\[
\nabla a(u,v)=\cases{ \displaystyle\frac{a(v)-a(u)}{v-u}, &\quad if $u\neq v$,
\vspace*{3pt}\cr
a'(u), &\quad if $u=v$} %
\]
for $u,v\in\R$ and consider the processes
\begin{eqnarray*}
\bigl(D^n_t\bigr)&=&\bigl(\nabla a\bigl(X^n_{\iota_n(t)},
X^n_t\bigr)\bigr),\qquad
\bigl(D^{n,n+1}_t \bigr)=\bigl(\nabla a\bigl(X^n_{t}, X^{n+1}_t
\bigr)\bigr),
\\
\bigl(A_t^n\bigr)&=&a\bigl(X^n_{\iota_n(t)}
\bigr). %
\end{eqnarray*}
In terms of the new notation, we have
%
\begin{eqnarray}
\label{eq3110-1} \dd U^{n,n+1}_{t}&=&D_{t-}^{n,n+1}
U^{n,n+1}_{t-} \,\dd Y_{t}+\varepsilon_n^{-1/2}
D^n_{t-} A_{t-}^n
(Y_{t-}-Y_{\iota
_{n}(t-)}) \,\dd Y_{t}
\nonumber\\[-8pt]\\[-8pt]\nonumber
&&{} -\varepsilon_n^{-1/2} D_{t-}^{n+1}
A_{t-}^{n+1} (Y_{t-}-Y_{\iota_{n+1}(t-)}) \,\dd
Y_{t}.
\end{eqnarray}
Clearly, the processes $(D_t^{n})$ and $(D_t^{n,n+1})$ converge in ucp
to $(D_t):=\break (a'(X_t))_{t\in[0,T]}$ and the processes $(A^n_t)$ to
$(A_t):=(a(X_t))_{t\in[0,T]}$. It often will be useful that the
processes $D^{n,n+1}$ and $D$ are uniformly bounded by the Lipschitz
constant of the coefficient $a$.

For technical reasons, we introduce a further approximation. For every
$\varepsilon>0$, we denote by $U^{n,n+1,\varepsilon
}=(U^{n,n+1,\varepsilon}_t)_{t\in[0,T]}$ the solution of the SDE
%
\begin{eqnarray}
\label{eq0602-3}
\dd U^{n,n+1,\varepsilon}_t &=& D_{t-}
U_{t-}^{n,n+1,\varepsilon} \,\dd Y_t
\nonumber\\[-8pt]\\[-8pt]\nonumber
&&{} + \varepsilon_n^{-1/2}
D_{t-} A_{t-} \sigma(W_{\iota
_{n+1}(t-)}-W_{\iota_{n}(t-)})
\,\dd Y^\varepsilon_{t}
\end{eqnarray}
with $U_0^{n,n+1,\varepsilon}=0$, where $Y^\varepsilon$ is as
in~(\ref{defYh}). Further, let $U^\varepsilon=(U^\varepsilon
_t)_{t\in[0,T]}$ denote the solution of
%
\begin{eqnarray}
\label{limitproceeps} U^\varepsilon_{t}&=&\int_{0}^{t}D_{s-}
U^\varepsilon_{s-} \,\dd Y_{s}+\sigma^2
\Upsilon\int_{0}^{t}D_{s-}A_{s-}
\,\dd B_{s}
\nonumber\\[-8pt]\\[-8pt]\nonumber
&&{} +\sum_{s\in(0,t]\dvtx  \Delta Y^\varepsilon_s\neq0}\sigma_s \mathcal
\xi_{s} D_{s-} A_{s-} \Delta Y_s^\varepsilon.
\end{eqnarray}

We\vspace*{1pt} will show that the processes $U^\varepsilon, U^{1,2,\varepsilon},
U^{2,3,\varepsilon}, \ldots$ are good approximations for the processes
$U, U^{1,2}, U^{2,3}, \ldots$ in the sense of Remark~\ref{re2901-1}.
As a consequence of Lemma~\ref{le2901-1}, we then get:

\begin{proposition}\label{prop1}
If for every $\varepsilon>0$,
\[
\bigl(Y,U^{n,n+1,\varepsilon}\bigr)\Rightarrow\bigl(Y, U^\varepsilon\bigr)\qquad
\mbox {in }\ID\bigl(\R^2\bigr), %
\]
then one has
\[
\bigl(Y,U^{n,n+1}\bigr)\Rightarrow(Y, U)\qquad\mbox{in }\ID\bigl(
\R^2\bigr). %
\]
\end{proposition}

The proof of the proposition is carried out in Section~\ref{sec21}.
It then remains to prove the following proposition which is the task of
Section~\ref{sec22}.

\begin{proposition}\label{prop2}For every $\varepsilon>0$,
\[
\bigl(Y, U^{n,n+1,\varepsilon}\bigr)\Rightarrow\bigl(Y, U^\varepsilon\bigr)\qquad
\mbox {in }\ID\bigl(\R^2\bigr). %
\]
\end{proposition}

\subsection{The approximations \texorpdfstring{$U^{n,n+1,\varepsilon}$}{$U^{n,n+1,varepsilon}$} are good}\label{sec21}
In this subsection, we prove Proposition~\ref{prop1}. By Lemma~\ref
{le2901-1}, it suffices to show that the approximations are good in the
sense of Remark~\ref{re2901-1}. In this section, we will work with an
additional auxiliary process: for $n\in\N$ and $\varepsilon>0$ we
denote by $\overline U^{n,n+1,\varepsilon}:=(\overline U^{n,n+1,\varepsilon
}_t)_{t\in[0,T]}$ the solution of
%
\begin{eqnarray}
\label{eq1902-2}
\qquad\dd\overline U^{n,n+1,\varepsilon}_{t}&=&D^{n,n+1}_{t-}
\overline U^{n,n+1,\varepsilon}_{t-} \,\dd Y_{t}+\varepsilon_n^{-1/2}
D^n_{t-} A^n_{t-}
\sigma(W_{t-}-W_{\iota_{n}(t-)}) \,\dd Y^{\varepsilon}_{t}
\nonumber\\[-8pt]\\[-8pt]\nonumber
&&{}-\varepsilon_n^{-1/2} D^{n+1}_{t-}
A^{n+1}_{t-} (W_{t-}-W_{\iota
_{n+1}(t-)}) \,\dd
Y^{\varepsilon}_{t}
\end{eqnarray}
with $\overline U^{n,n+1,\varepsilon}_0=0$.

\begin{lemma}\label{lemmaappgood}
For every $\delta,\varepsilon>0$, we have:
\begin{longlist}
\item[1.] $\lim_{\varepsilon\downarrow0}\limsup_{n\to\infty}\E
[\sup_{t\in[0,T]}\llvert  U^{n,n+1}_t-\overline U^{n,n+1,\varepsilon}_t\rrvert  ^2 ]=0$,\vspace*{2pt}

\item[2.]
$\lim_{n\to\infty}\P (\sup_{t\in[0,T]}\llvert   \overline U^{n,n+1,\varepsilon}_{t}- U^{n,n+1,\varepsilon}_{t}\rrvert   >\delta
 )=0$,\vspace*{2pt}

\item[3.]
$\lim_{\varepsilon\downarrow0}\P (\sup_{t\in
[0,T]}\llvert  U_t-U^{\varepsilon}_t\rrvert  >\delta )=0$.
\end{longlist}
\end{lemma}

It is straightforward to verify that Lemma~\ref{lemmaappgood} implies
that the approximations are good.

\begin{pf*}{Proof of Lemma \ref{lemmaappgood}}
(1)
Recalling~(\ref{eq3110-1}) and~(\ref{eq1902-2}) and noting that
$D^{n,n+1}$ is uniformly bounded, we conclude with Lemma~\ref
{estimatesdes2} that the first statement is true if
%
\begin{eqnarray}\label{goodeq2}
&&\lim_{\varepsilon\downarrow0}\limsup_{n\to\infty}
\varepsilon _n^{-1} \E \biggl[\sup_{t\in[0,T]}
\biggl\llvert \int_0^tD^n_{s-}A^n_{s-}(Y_{s-}-Y_{\iota_n(s-)})
\,\dd Y_s
\nonumber\\[-8pt]\\[-8pt]\nonumber
&&\hspace*{106pt}{} -\sigma\int_0^tD^n_{s-}A^n_{s-}(W_{s-}-W_{\iota_n(s-)})
\,\dd Y^{\varepsilon}_s\biggr\rrvert ^2 \biggr]=0.
\end{eqnarray}
Let $M^{\varepsilon}$ denote the martingale $Y-Y^{\varepsilon}$. The
above term can be estimated against the sum of
%
\begin{equation}
\label{eq1202-1} \qquad\varepsilon_n^{-1}\E \biggl[\sup
_{t\in[0,T]}\biggl\llvert \int_0^tD^n_{s-}A^n_{s-}(Y_{s-}-Y_{\iota_n(s-)}-
\sigma W_{s-}+\sigma W_{\iota_n(s-)}) \,\dd Y_s\biggr
\rrvert ^2 \biggr]
\end{equation}
and
%
\begin{equation}
\label{eq1202-2} \varepsilon_n^{-1}\sigma^2\E
\biggl[\sup_{t\in[0,T]}\biggl\llvert \int_0^tD^n_{s-}A^n_{s-}(W_{s-}-W_{\iota_n(s-)})
\,\dd M^{\varepsilon
}_s\biggr\rrvert ^2 \biggr].
\end{equation}
We start with estimating the former expression. For $t\in[0,T]$, one has
\begin{eqnarray*}
Y_{t}-Y_{\iota_{n}(t)} &=& \sigma(W_{t}-W_{\iota
_{n}(t)})+M^{h_{n}}_{t}-M^{h_{n}}_{\iota_{n}(t)}
\\
&&{} +\biggl(b-\int_{B(0,h_{n})^c}x \nu(\dd x) \biggr) \bigl(t-
\iota_n(t)\bigr).
\end{eqnarray*}
By Lemma~\ref{driftdisappear}, one has
\begin{eqnarray*}
&& \varepsilon_n^{-1} \E \bigl[\llvert Y_{t}-Y_{\iota_n(t)}-
\sigma W_{t}+\sigma W_{\iota_{n}(t)}\rrvert ^2\mid
\iota_n \bigr]
\\
&&\qquad \leq  2\int_{B(0,h_n)}x^2
\nu(\dd x) +2 \biggl(b-\int_{B(0,h_{n})^c}x \nu(\dd x)
\biggr)^2\varepsilon_n
\\
&&\qquad =:\delta_n\to0
\end{eqnarray*}
as $n\to\infty$.
Further, by Lemma~\ref{estimateintegral} and the uniform boundedness
of $D^n$, there is a constant $\kappa_1$ not depending on $n$ such that
%
\begin{eqnarray}
\label{eq1202-3}
\qquad&& \varepsilon_n^{-1} \E \biggl[\sup
_{t\in[0,T]}\biggl\llvert \int_0^tD^n_{s-}A^n_{s-}(Y_{s-}-Y_{\iota_n(s-)}-
\sigma W_{s-}+\sigma W_{\iota_n(s-)}) \,\dd Y_s\biggr
\rrvert ^2 \biggr]\nonumber
\\
&&\qquad \leq\kappa_1 \varepsilon_n^{-1} \int
_0^T\E\bigl[\bigl(A^n_{s-}
\bigr)^2 (Y_{s-}-Y_{\iota_n(s-)}-\sigma
W_{s-}+\sigma W_{\iota_{n}(s-)})^2\bigr] \,\dd s
\\
&&\qquad \leq\kappa_1 \delta_n \int_0^T
\E\bigl[\bigl\llvert A^n_{s-}\bigr\rrvert ^2
\bigr] \,\dd s,\nonumber
\end{eqnarray}
where we have used conditional independence of $A_{s-}^n$ and
$Y_{s-}-Y_{\iota_n(s-)}-\sigma W_{s-}+\sigma W_{\iota_{n}(s-)}$ given
$\iota_n$ in the last transformation.
By Lemma~\ref{boundapp} and the Lipschitz continuity of $a$, the
latter integral is uniformly bounded over all $n\in\N$ so that~(\ref
{eq1202-1}) tends to zero as $n\to\infty$.

Next, consider~(\ref{eq1202-2}). Note that $M^\varepsilon$ is a L\'
evy martingale with triplet $(0,0, \nu\mid _{B(0,\varepsilon)})$. By
Lemma~\ref{estimateintegral} and the uniform boundedness of $D^n$,
there exists a constant $\kappa_2$ not depending on $\varepsilon$ and
$n$ such that
%
\begin{eqnarray}
\label{eq1902-4}
&& \varepsilon_n^{-1} \E \biggl[\sup
_{t\in[0,T]}\biggl\llvert \int_0^tD^n_{s-}A^n_{s-}(W_{s-}-W_{\iota_n(s-)})
\,\dd M^{\varepsilon
}_s\biggr\rrvert ^2 \biggr]\nonumber
\\
&&\qquad \leq\kappa_2 \varepsilon_n^{-1} \int
_{B(0,\varepsilon)}x^2 \nu (\dd x)\int_0^T
\E\bigl[\bigl\llvert A^n_{s-}\bigr\rrvert ^2
\llvert W_{s-}-W_{\iota_n(s-)}\rrvert ^2\bigr] \,\dd s
\\
&&\qquad \leq\kappa_2 \int_{B(0,\varepsilon)}x^2 \nu(\dd
x)\int_0^T\E \bigl[\bigl\llvert
A^n_{s-}\bigr\rrvert ^2\bigr] \,\dd s,\nonumber
\end{eqnarray}
where we used in the last step that conditionally on $\iota_n$ the
random variables $A^n_{s-}$ and $W_{s-}-W_{\iota_n(s-)}$ are
independent and $\E[(W_{s-}-W_{\iota_n(s-)})^2\mid \iota_n]=s-\iota
_n(s)\leq\varepsilon_n$.
As noted above, $\int_0^T\E[\llvert  A^n_{s-}\rrvert  ^2] \,\dd s$ is uniformly
bounded, and hence~(\ref{eq1202-2}) tends uniformly to zero over all
$n\in\N$ as $\varepsilon\downarrow0$.

(2)
We will use Lemma~\ref{estimationsde} to prove that
%
\begin{equation}
\label{goodeq3} \overline U^{n,n+1,\varepsilon}-U^{n,n+1,\varepsilon}\to0\qquad\mbox{in ucp, as
}n\to\infty.
\end{equation}
We rewrite the SDE~(\ref{eq0602-3}) as
\begin{eqnarray*}
\dd U^{n,n+1,\varepsilon}_t&=&D_{t-} U_{t-}^{n,n+1,\varepsilon}
\,\dd Y_t+ \varepsilon_n^{-1/2} D_{t-}
A_{t-} \sigma(W_{t-}-W_{\iota
_{n}(t-)}) \,\dd
Y^\varepsilon_{t}
\\
&&{}-\varepsilon_n^{-1/2} D_{t-} A_{t-}
\sigma(W_{t-}-W_{\iota
_{n+1}(t-)}) \,\dd Y^\varepsilon_{t}.
\end{eqnarray*}
Recalling~(\ref{eq1902-2}), it suffices by part one of Lemma~\ref
{estimationsde} to show that:
\begin{longlist}[2.]
\item[1.] $D^{n,n+1}\to D$, in ucp,

\item[2.] $\varepsilon_n^{-1/2}\int_0^{\cdot
}(D^n_{s-}A^n_{s-}-D_{s-}A_{s-})(W_{s-}-W_{\iota_n(s-)}) \,\dd
Y^{\varepsilon}_s\to0$, in ucp,

\item[3.]  the families $(\sup_{t\in[0,T]}\llvert  D_t^{n,n+1}\rrvert\dvtx n\in\N)$ and
\[
\biggl(\varepsilon_n^{-1/2}\sup_{t\in[0,T]}
\biggl\llvert \int_0^{t}D^n_{s-}A^n_{s-}(W_{s-}-W_{\iota_n(s-)})
\,\dd Y^{\varepsilon
}_s\biggr\rrvert\dvtx n\in\N \biggr) %
\]
are tight.
\end{longlist}
The\vspace*{1pt} tightness of $(\sup_{t\in[0,T]}\llvert  D_t^{n,n+1}\rrvert\dvtx n\in\N)$ follows
by uniform boundedness. Further, the tightness of the second family
follows by observing that in analogy to the proof of~(1) one has
%
\begin{eqnarray*}
&& \varepsilon_n^{-1} \E \biggl[\sup_{t\in[0,T]}
\biggl\llvert \int_0^tD^n_{s-}A^n_{s-}(W_{s-}-W_{\iota_n(s-)})
\,\dd Y^{\varepsilon
}_s\biggr\rrvert ^2 \biggr]
\\
&&\qquad \leq\kappa_3 \varepsilon_n^{-1} \int
_0^T\E \bigl[\bigl\llvert A^n_{s-}
\bigr\rrvert ^2\llvert W_{s-}-W_{\iota_n(s-)}\rrvert
^2\bigr] \,\dd s\leq\kappa_3 \int_0^T
\E\bigl[\bigl\llvert A^n_{s-}\bigr\rrvert ^2
\bigr] \,\dd s
\end{eqnarray*}
for an appropriate constant $\kappa_3$ not depending on $n$.
Furthermore, convergence $D^{n,n+1}\to D$ follows from ucp convergence
of $ X^n\to X$ and Lipschitz continuity of $a$.
To show the remaining property, we let $\delta>0$ and $T_{n,\delta}$
denote the stopping time
\[
T_{n,\delta}=\inf\bigl\{s\in[0,T]\dvtx  \bigl\llvert D_s^nA^n_s-D_sA_s
\bigr\rrvert \geq\delta\bigr\}. %
\]
Then by Lemma~\ref{estimateintegral}, there exists a constant $\kappa
_4$ not depending on $n$ and $\delta$ with
\begin{eqnarray*}
&& \E \biggl[\sup_{t\in[0,T\wedge T_{n,\delta}]} \varepsilon_n^{-1}
\biggl(\int_0^t \bigl(D^n_{s-}A^n_{s-}-D_{s-}A_{s-}
\bigr) (W_{s-}-W_{\iota
_n(s-)}) \,\dd Y^{\varepsilon}_s
\biggr)^2 \biggr]
\\
&&\qquad \leq\kappa_4 \delta^2 \varepsilon_n^{-1}
\int_0^T \E \bigl[(W_{s-}-W_{\iota_n(s-)})^2
\bigr] \,\dd s\leq\kappa_4 \delta^2 T.
\end{eqnarray*}
Since for any $\delta>0$, $\P(T_{n,\delta}=\infty)\to1$ by ucp
convergence $D^nA^n-DA\to0$, we immediately get the remaining property
by choosing $\delta>0$ arbitrarily small and applying the Markov
inequality.

(3)
The proof of the third statement can be achieved by a simplified
version of the proof of the first statement. It is therefore omitted.
\end{pf*}

\subsection{Weak convergence of \texorpdfstring{$U^{n,n+1,\varepsilon}$}{$U^{n,n+1,varepsilon}$}}\label{sec22}

In this subsection, we prove Proposition~\ref{prop2} for fixed
$\varepsilon>0$. We first outline the proof.
We will make use of results of~\cite{jacodprotter} summarised in the
\hyperref[append]{Appendix}; see Section~\ref{sec1803-01}. We consider processes
$Z^{n,\varepsilon}=(Z^{n,\varepsilon}_t)_{t\in[0,T]}$ and
$Z^{\varepsilon}=(Z^{\varepsilon}_t)_{t\in[0,T]}$ given by
%
\begin{equation}
\label{eq3101-2} Z^{n,\varepsilon}_t =\varepsilon_n^{-1/2}
\int_0^t (W_{\iota_{n+1}
(s-)}-W_{\iota_{n} (s-)})
\,\dd Y^\varepsilon_s
\end{equation}
and
%
\begin{equation}
\label{eq3101-3} Z^\varepsilon_t= \Upsilon B_{t} +\sum
_{s\in(0,t]\dvtx  \llvert  \Delta Y_s\rrvert  \geq\varepsilon}\frac{\sigma
_s}{\sigma} \mathcal\xi_{s}
\Delta Y_{s},
\end{equation}
where $(\sigma_s)$ and $(\xi_s)$ are the marks of the point
process~$\Pi$ as introduced in Section~\ref{sec1803-02}.


In view of Theorem~\ref{weakproperty2}, the statement of
Proposition~\ref{prop2} follows, if we show that
%
\begin{eqnarray}
\biggl(Y,\int_0^\cdot D_{t-} \,\dd
Y_t, \int_0^\cdot D_{t-}
A_{t-} \,\dd Z_t^{n,\varepsilon} \biggr) \Rightarrow
\biggl(Y, \int_0^\cdot D_{t-} \,\dd
Y_t, \int_0^\cdot D_{t-}
A_{t-} \,\dd Z_t^{\varepsilon
} \biggr)\nonumber
\\
\eqntext{\mbox{in }\ID \bigl(\R^3\bigr).}
\end{eqnarray}
Further, by Theorem~\ref{weakproperty1}, this statement follows once
we showed that\break $(Z^{n,\varepsilon}\dvtx n\in\N)$ is uniformly tight and
%
\begin{equation}
\label{eq1703-2} \bigl( Y,D, DA, Z^{n,\varepsilon} \bigr) \Rightarrow \bigl( Y,
D,DA, Z^{\varepsilon} \bigr)\qquad\mbox{in }\ID\bigl(\R^4\bigr).
\end{equation}

We first prove that $((Y,D,DA,Z^{n,\varepsilon})\dvtx n\in\N)$ is tight
which shows that, in particular, $(Z^{n,\varepsilon}\dvtx n\in\N)$ is
uniformly tight; see Lemma~\ref{proptightness}. Note that $(Y,D,DA)$
is $\sigma(Y)$-measurable. To identify the limit and complete the
proof of~(\ref{eq1703-2}), it suffices to prove stable convergence
\[
Z^{n,\varepsilon}\stackrel{\operatorname{stably}} {\Longrightarrow
}Z^\varepsilon %
\]
with respect to the $\sigma$-field $\sigma(Y)$; see Section~\ref
{sec1803-01} in the \hyperref[append]{Appendix} for a brief introduction of stable
convergence. The latter statement is equivalent to
\[
\bigl(Y,Z^{n,\varepsilon}\bigr)\Rightarrow\bigl(Y,Z^\varepsilon\bigr)\qquad
\mbox{in } \ID(\R)\times\ID(\R), %
\]
by Theorem~\ref{stablethm1}. We prove the stronger statement that
this is even true in the finer topology~$\ID(\R^2)$: the sequence
$((Y,Z^{n,\varepsilon})\dvtx n\in\N)$ is tight by Lemma~\ref
{proptightness} and we will prove convergence of finite-dimensional
marginals in Lemma~\ref{le2901-2}. The proof of the latter lemma is
based on a perturbation result provided by Lemma~\ref{le3101-2}.

\begin{lemma}\label{proptightness}For $\varepsilon>0$, the family
$((Y,D,DA,Z^{n,\varepsilon})\dvtx n\in\N)$ taking values in $\ID(\R
^4)$ is tight. In particular, $(Z^{n,\varepsilon}\dvtx n\in\N)$ is
uniformly tight.
\end{lemma}

\begin{pf}
One has by Lemma~\ref{estimateintegral}
\[
\E \Bigl[\sup_{t\in[0,T]} \bigl(Z^{n,\varepsilon}_t
\bigr)^2 \Bigr] \leq \kappa_1 \varepsilon_n^{-1}
\int_0^t \E\bigl[(W_{\iota
_{n+1}(t-)}-W_{\iota_{n}(t-)})^2
\bigr] \,\dd t\leq\kappa_1 %
\]
for an appropriate constant $\kappa_1$
so that by the Markov inequality
\[
\lim_{K\to\infty} \sup_{n\in\N} \P \Bigl(\sup
_{t\in[0,T]} \llvert Y_t\rrvert \vee\bigl\llvert
Z^{n,\varepsilon}_t\bigr\rrvert \vee\llvert D_t\rrvert
\vee\llvert D_tA_t\rrvert \geq K \Bigr)=0. %
\]
It remains to verify Aldous' criterion for tightness \cite{jacodshibook}, Theorem~VI.4.5, which can be checked componentwise. It
is certainly fulfilled for $Y$, $A$ and $DA$ and it remains to show
that for every $K>0$ there exists for every $\delta>0$ a constant
$c_\delta>0$ such that for arbitrary stopping times $S_1,S_2,\ldots$
\[
\limsup_{n\to\infty} \P \Bigl(\sup_{t\in[S_n,(S_n+\delta)\wedge
T]} \bigl
\llvert Z^{n,\varepsilon}_t-Z^{n,\varepsilon}_{S_n}\bigr
\rrvert \geq K \Bigr)\leq c_\delta %
\]
and $\lim_{\delta\downarrow0} c_\delta=0$.

First, suppose that $S_1,S_2,\ldots$ denote stopping times taking
values in the respective sets~$ \varepsilon_n\Z$. Then as above
%
\begin{eqnarray}\label{eq2602-1}
&& \E \Bigl[\sup_{t\in[S_n,(S_n+\delta)\wedge T]} \bigl\llvert
Z^{n,\varepsilon
}_t-Z^{n,\varepsilon}_{S_n}\bigr\rrvert
^2 \Bigr]\nonumber
\\
&&\qquad  \leq \kappa_1 \varepsilon _n^{-1}
\int_0^T \E \bigl[\1_{[S_n,(S_n+\delta)]}(t)
(W_{\iota
_{n+1}(t-)}-W_{\iota_{n}(t-)})^2 \bigr] \,\dd t\nonumber
\\
&&\qquad \leq \kappa_1 \varepsilon_n^{-1} \int
_0^T\E \bigl[ \1 _{[S_n,(S_n+\delta)]}\bigl(
\iota_{n}(t)\bigr) (W_{\iota_{n+1}(t-)}-W_{\iota
_{n}(t-)})^2\bigr] \,\dd t
\\
&&\qquad \leq \kappa_1 \E \biggl[\int_0^T
\1_{[S_n,(S_n+\delta)]}\bigl(\iota _{n}(t)\bigr) \,\dd t \biggr]\nonumber
\\
&&\qquad \leq
\kappa_1 (\varepsilon_n+\delta)\to\kappa _1
\delta,\nonumber
\end{eqnarray}
where we have used that $\E[(W_{\iota_{n+1}(t-)}-W_{\iota
_{n}(t-)})^2\mid \cF_{\iota_n(t)}]\leq\varepsilon_n$ and $\1
_{[S_n,(S_n+\delta)]}(t)$ is $\cF_{\iota_n(t)}$-measurable. It
remains to estimate for general stopping times $S_1,S_2,\ldots$
\[
\E \Bigl[\sup_{t\in[S_n,\overline S_n]} \bigl\llvert Z^{n,\varepsilon
}_t-Z^{n,\varepsilon}_{S_n}
\bigr\rrvert ^2 \Bigr], %
\]
where $\overline S_n=\inf[S_n,\infty)\cap\varepsilon_n\Z$.
As in~(\ref{eq2602-1}), we conclude with $\overline S_n-S_n\leq\varepsilon
$ that
\begin{eqnarray*}
&& \E \Bigl[\sup_{t\in[S_n,\overline S_n]} \bigl\llvert Z^{n,\varepsilon
}_t-Z^{n,\varepsilon}_{S_n}
\bigr\rrvert ^2 \Bigr]
\\
&&\qquad  \leq\kappa_1 \varepsilon
_n^{-1} \E \biggl[ \int_0^T
\1_{[S_n,\overline S_n]}(t) (W_{\iota
_{n+1}(t-)}-W_{\iota_{n}(t-)})^2 \,\dd t
\biggr]
\\
&&\qquad \leq\kappa_1 \E \Bigl[ \mathop{\sup_{k=1,\ldots,\varepsilon
_n^{-1}}}_{s,t\in[(k-1)\varepsilon_n^{-1},k\varepsilon_n^{-1})}
\llvert W_{s}-W_{t}\rrvert ^2 \Bigr]\to0.
\end{eqnarray*}
By the Markov inequality, this estimate together with~(\ref{eq2602-1})
imply Aldous' criterion.
\end{pf}

To control perturbations, we will use the following lemma.

\begin{lemma}\label{le3101-2}For $j=1,2$, let $(\alpha^{(j)}_t)_{t\in
[0,T]}$ and $(\beta^{(j)}_t)_{t\in[0,t]}$ optional processes being
square integrable with respect to $\P\otimes\ell_{[0,T]}$ and let
\[
\Upsilon^{n,j}_t= \varepsilon_n^{-1/2}
\int_0^t \bigl(\overline W^{(j)}_{\iota
_{n+1}(s-)}-
\overline W^{(j)}_{\iota_{n}(s-)}\bigr) \,\dd\overline Y^{(j)}_s,
\]
where
\[
\overline W^{(j)}_t=W_t+\int_0^t
\alpha^{(j)}_s \,\dd s,\qquad \overline Y^{(j)}_t=
M_t+\int_0^t \beta^{(j)}_s
\,\dd s
\]
and
\[
M_t= \sigma W_t+\int
_{(0,t]\times B(0,\varepsilon)^c} x \,\dd\overline\Pi(s,x). %
\]
For $t\in\cD=\bigcup_{n\in\N} \varepsilon_n\Z\cap[0,T]$, the
sequences $(\Upsilon^{n,1}_t)_{n\in\N}$ and $(\Upsilon
_t^{n,2})_{n\in\N}$ are equivalent in probability, that is, for every
$\delta>0$
\[
\lim_{n\to\infty} \P \bigl( \bigl\llvert \Upsilon^{n,1}_t-
\Upsilon ^{n,2}_t\bigr\rrvert >\delta \bigr)=0. %
\]
\end{lemma}

\begin{pf}
We prove the statement in three steps.

\begin{longlist}[1\textit{st step}.]
\item[1\textit{st step}.] First, we show a weaker perturbation estimate. Using
the bilinearity of the stochastic integral, we get that
%
\begin{eqnarray}
\label{eq0302-1} \Upsilon^{n,1}_t-\Upsilon^{n,2}_t
&=& \varepsilon_n^{-1/2} \int_0^t
\int_{\iota_n(s-)}^{\iota_{n+1}(s-)} \bigl(\alpha_u^{(1)}-
\alpha _u^{(2)}\bigr) \,\dd u \,\dd M_s\nonumber
\\
&&{} +\varepsilon_n^{-1/2}\int_0^t
(W_{\iota_{n+1}(s-)}-W_{\iota
_n(s-)}) \bigl(\beta^{(1)}_s-
\beta^{(2)}_s\bigr) \,\dd s
\nonumber\\[-8pt]\\[-8pt]\nonumber
&&{} + \varepsilon_n^{-1/2} \int_0^t
\int_{\iota_n(s-)}^{\iota
_{n+1}(s-)} \bigl(\alpha_u^{(1)}-
\alpha_u^{(2)}\bigr) \,\dd u\, \beta^{(1)}_s
\,\dd s
\\
&&{} + \varepsilon_n^{-1/2} \int_0^t
\int_{\iota_n(s-)}^{\iota
_{n+1}(s-)} \alpha_u^{(2)}
\,\dd u \bigl(\beta^{(1)}_s-\beta^{(2)}_s
\bigr) \,\dd s.\nonumber
\end{eqnarray}

We analyse the terms individually. By It\^o's isometry, the fact that
$s-\varepsilon_n\leq\iota_n(s-)\leq\iota_{n+1}(s-)\leq s$ and
Fubini's theorem one has that for $\kappa=\sigma^2+\int_{B(0,\varepsilon)^c} x^2 \nu(\dd x)$
%
\begin{eqnarray}
\label{eq1501-1}
&& \E \biggl[  \biggl(\varepsilon_n^{-1/2} \int
_0^t \int_{\iota
_n(s-)}^{\iota_{n+1}(s-)}
\bigl(\alpha_u^{(1)}-\alpha_u^{(2)}
\bigr) \,\dd u \,\dd M_s \biggr)^2 \biggr]\nonumber
\\
&&\qquad = \kappa \varepsilon_n^{-1} \E \biggl[\int
_0^t \biggl(\int_{\iota_n(s-)}^{\iota_{n+1}(s-)}
\bigl(\alpha_u^{(1)}-\alpha_u^{(2)}
\bigr) \,\dd u \biggr)^2 \,\dd s \biggr]\nonumber
\\
&&\qquad \leq\kappa \E \biggl[\int_0^t \int
_{\iota_n(s-)}^{\iota
_{n+1}(s-)} \bigl(\alpha_u^{(1)}-
\alpha_u^{(2)}\bigr)^2 \,\dd u \,\dd s \biggr]
\\
&&\qquad \leq\kappa \E \biggl[\int_0^t \int
_{(s-\varepsilon_n)\vee0}^{s} \bigl(\alpha_u^{(1)}-
\alpha_u^{(2)}\bigr)^2 \,\dd u \,\dd s \biggr]\nonumber
\\
&&\qquad \leq\kappa \varepsilon_n \E \biggl[\int_0^t
\bigl(\alpha _s^{(1)}-\alpha_s^{(2)}
\bigr)^2 \,\dd s \biggr].\nonumber
\end{eqnarray}
By the Cauchy--Schwarz inequality and Fubini, it follows that the
second term satisfies
\begin{eqnarray*}
&& \E \biggl[ \varepsilon_n^{-1/2}\biggl\llvert \int
_0^t (W_{\iota
_{n+1}(s-)}-W_{\iota_n(s-)})
\bigl(\beta^{(1)}_s-\beta^{(2)}_s
\bigr) \,\dd s\biggr\rrvert \biggr]
\\
&&\qquad \leq\varepsilon_n^{-1/2} \E \biggl[\int
_0^t (W_{\iota
_{n+1}(s-)}-W_{\iota_n(s-)})^2
\,\dd s \biggr]^{1/2} \E \biggl[\int_0^t
\bigl(\beta^{(1)}_s-\beta^{(2)}_s
\bigr)^2 \,\dd s \biggr]^{1/2}
\\
&&\qquad \leq t \E \biggl[\int_0^t \bigl(
\beta^{(1)}_s-\beta^{(2)}_s
\bigr)^2 \,\dd s \biggr]^{1/2},
\end{eqnarray*}
where we have used in the last step that $\iota^{n+1}_{s-}-\iota
^n_{s-}$ is independent of the Brownian motion and smaller or equal to
$\varepsilon_n$. The third term is estimated similarly as the first term:
\begin{eqnarray*}
&& \E \biggl[ \varepsilon_n^{-1/2}\biggl\llvert \int
_0^t \int_{\iota
_n(s-)}^{\iota_{n+1}(s-)}
\bigl(\alpha_u^{(1)}-\alpha_u^{(2)}
\bigr) \,\dd u\, \beta^{(1)}_s \,\dd s\biggr\rrvert \biggr]
\\
&&\qquad \leq\varepsilon_n^{1/2} \E \biggl[\int
_0^t \biggl(\int_{\iota
_n(s-)}^{\iota_{n+1}(s-)}
\bigl(\alpha_u^{(1)}-\alpha_u^{(2)}
\bigr) \,\dd u \biggr)^2 \,\dd s \biggr]^{1/2} \E \biggl[\int
_0^T \bigl(\beta ^{(1)}_s
\bigr)^2 \,\dd s \biggr]^{1/2}
\\
&&\qquad \leq\E \biggl[\int_0^t\int
_{\iota_n(s-)}^{\iota_{n+1}(s-)} \bigl(\alpha_u^{(1)}-
\alpha_u^{(2)}\bigr)^2 \,\dd u \,\dd s
\biggr]^{1/2} \E \biggl[\int_0^T
\bigl(\beta^{(1)}_s\bigr)^2 \,\dd s
\biggr]^{1/2}
\\
&&\qquad \leq\varepsilon_n^{1/2} \E \biggl[\int
_0^t \bigl(\alpha _s^{(1)}-
\alpha_s^{(2)}\bigr)^2 \,\dd s
\biggr]^{1/2} \E \biggl[\int_0^T
\bigl(\beta^{(1)}_s\bigr)^2 \,\dd s
\biggr]^{1/2}.
\end{eqnarray*}
In complete analogy, the fourth term satisfies
\begin{eqnarray*}
&& \E \biggl[  \varepsilon_n^{-1/2}\biggl\llvert \int
_0^t \int_{\iota
_n(s-)}^{\iota_{n+1}(s-)}
\alpha_u^{(2)} \,\dd u \bigl(\beta ^{(1)}_s-
\beta^{(2)}_s\bigr) \,\dd s\biggr\rrvert \biggr]
\\
&&\qquad \leq\varepsilon_n^{1/2} \E \biggl[\int
_0^t\bigl(\alpha_s^{(2)}
\bigr)^2 \,\dd s \biggr]^{1/2} \E \biggl[\int
_0^T \bigl(\beta^{(1)}_s-
\beta ^{(2)}_s\bigr)^2 \,\dd s
\biggr]^{1/2}.
\end{eqnarray*}
By the Markov inequality, the first, third and fourth term of~(\ref
{eq0302-1}) tend to zero in probability as $n\to\infty$.

\item[2\textit{nd step}.] Next, we analyse the case where $\beta^{(2)}=0$ and
$\beta:=\beta^{(1)}$ is simple in the following sense. There exist
$l\in\N$, increasingly ordered times $0=t_0,t_1,\ldots,t_l=t\in\cD
=\bigcup_{n\in\N} \varepsilon_n\mathbb{Z}\cap[0,T]$
such that $\beta$ is almost surely constant on each of the time
intervals $[t_0,t_1),\ldots,[t_{l-1},t_l)$.
For $n\in\N$ and $j=1,\ldots,l$, we let
\[
M_{j,n}:=\varepsilon_n^{-1/2} \int
_{t_{j-1}}^{t_j} (W_{\iota
_{n+1}(s-)}-W_{\iota_{n}(s-)}) \,\dd
s. %
\]
We suppose that $n\in\N$ is sufficiently large to ensure that $\{
t_1,\ldots,t_l\}\subset\varepsilon_n\Z$. The Brownian motion $W$ is
independent of $\Pi$ so that for $u,s\in[0,t]$
\begin{eqnarray*}
&& \E\bigl[(W_{\iota_{n+1}(s-)}-W_{\iota_{n}(s-)}) (W_{\iota
_{n+1}(u-)}-W_{\iota_{n}(u-)})
\mid \Pi\bigr]
\\
&&\qquad =\ell\bigl(\bigl[\iota_n(s-), \iota_{n+1}(s-)\bigr]\cap
\bigl[\iota_n(u-), \iota _{n+1}(u-)\bigr]\bigr)
\\
&&\qquad \leq\varepsilon_n \1_{\{\llvert  s-u\rrvert  \leq\varepsilon_n\}}.
\end{eqnarray*}
Consequently, we obtain with Fubini that
\begin{eqnarray*}
\E\bigl[M_{j,n}^2\bigr] &=&\varepsilon_n^{-1}
\E \biggl[\int_{t_{j-1}}^{t_j}\int_{t_{j-1}}^{t_j}
(W_{\iota_{n+1}(s-)}-W_{\iota
_{n}(s-)}) (W_{\iota_{n+1}(u-)}-W_{\iota_{n}(u-)}) \,\dd s
\,\dd u \biggr]
\\
&\leq& 2 \varepsilon_n (t_j-t_{j-1}).
\end{eqnarray*}
Since $M_{j,n}$ is independent of $\cF_{t_{j-1}}$ and has mean zero,
we conclude that $(\sum_{j=1}^k \beta_{t_{j-1}} M_{j,n})_{k=0,\ldots,l}$ is a square integrable martingale so that
\begin{eqnarray*}
&& \E \biggl[ \biggl(\varepsilon_n^{-1/2}\int
_0^t (W_{\iota
_{n+1}(s-)}-W_{\iota_n(s-)})
\beta_s \,\dd s \biggr)^2 \biggr]
\\
&&\qquad =\sum
_{j=1}^l \E\bigl[\beta_{t_{j-1}}^2
M_{j,n}^2\bigr]
=\sum_{j=1}^l \E\bigl[
\beta_{t_{j-1}}^2\bigr] \E\bigl[M_{j,n}^2
\bigr]
\\
&&\qquad \leq 2\varepsilon_n \sum_{j=1}^l
\beta_{t_{j-1}}^2 (t_j-t_{j-1})
=2\varepsilon_n \E \biggl[\int_0^t
\beta_s^2 \,\dd s \biggr].
\end{eqnarray*}

\item[3\textit{rd step}.] We combine the first and second step. Let $\alpha
^{(2)}$ and $\beta^{(2)}$ be as in the statement of the theorem and
let $\delta>0$ be arbitrary. The simple functions as defined in step
two are dense in the space of previsible processes with finite
$L^2$-norm with respect\vspace*{1pt} to $\P\otimes\ell_{[0,T]}$. By part one, we
can choose $\alpha^{(1)}=0$ and a simple process $\beta^{(1)}$ such that
\[
\P\bigl( \bigl\llvert \Upsilon^{n,2}_t-
\Upsilon^{n,1}_t\bigr\rrvert \geq\delta/2\bigr)\leq\delta/2
\]
for $n$ sufficiently large.
Next, let $\Upsilon^{n,0}$ denote the process that is obtained in
analogy to $\Upsilon^{n,1}$ and $\Upsilon^{n,2}$ when choosing
$\alpha=\beta=0$. By the second step, $(\Upsilon^{n,1}_t\dvtx n\in\N)$
and $(\Upsilon^{n,0}_t\dvtx n\in\N)$ are asymptotically equivalent in
probability implying that
\[
\P\bigl(\bigl\llvert \Upsilon^{n,1}_t-\Upsilon^{n,0}_t
\bigr\rrvert \geq\delta/2\bigr)\leq\delta/2 %
\]
for sufficiently large $n\in\N$. Altogether, we arrive at
\[
\P\bigl( \bigl\llvert \Upsilon^{n,2}_t-
\Upsilon^{n,0}_t\bigr\rrvert \geq\delta\bigr)\leq\delta
\]
for sufficiently large $n\in\N$. Since $\delta>0$ is arbitrary,
$(\Upsilon^{n,2}_t\dvtx n\in\N)$ and $(\Upsilon^{n,0}_t\dvtx n\in\N)$ are
equivalent in probability. The general statement follows by
transitivity of equivalence in probability.\quad\qed
\end{longlist}\noqed
\end{pf}

\begin{lemma}\label{le2901-2} For any finite subset $\IT\subset\cD=
\bigcup_{n\in\N} \varepsilon_{n}  \N_0$, one
has convergence
\[
\bigl(Y_t,Z^{n,\varepsilon}_t\bigr)_{t\in\IT}
\Rightarrow\bigl(Y_t,Z^\varepsilon _t
\bigr)_{t\in\IT}.
\]
\end{lemma}

\begin{pf}
1\textit{st step}. In the first step, we derive a simpler
sufficient criterion which implies the statement. Fix $l\in\N$,
increasing times $0=t_0\leq t_1<\cdots<t_l\leq T$ and consider $\IT=\{
t_1,\ldots,t_l\}$.
The statement follows if for $A\in\sigma(Y_t\dvtx t\in\IT)$ and
continuous compactly supported $f\dvtx \R^l\to\R$
\[
\E\bigl[\1_A f\bigl(Z_{t_1}^{n,\varepsilon},
\ldots,Z_{t_l}^{n,\varepsilon
}\bigr)\bigr]\to\E\bigl[\1_A f
\bigl(Z_{t_1}^{\varepsilon},\ldots,Z_{t_l}^{\varepsilon}
\bigr)\bigr]. %
\]
By the Stone--Weierstrass theorem, the linear hull of functions of the form
\[
\R^l\to\R,\qquad x\mapsto f_1(x_1)\times \cdots \times f_l(x_l) %
\]
with continuous compactly supported functions $f_1,\ldots,f_l\dvtx \R\to\R
$ is dense in the space of compactly supported continuous functions on
$\R^l$ equipped with supremum norm. Hence, it suffices to verify that
%
\begin{eqnarray}
\label{eq1401-1} && \E\bigl[\1_A f_1\bigl(Z_{t_1}^{n,\varepsilon}
\bigr)\cdots f_l\bigl(Z_{t_l}^{n,\varepsilon
}-Z_{t_{l-1}}^{n,\varepsilon}
\bigr)\bigr]
\nonumber\\[-8pt]\\[-8pt]\nonumber
&&\qquad \to\E\bigl[\1_A f_1\bigl(Z_{t_1}^{\varepsilon
}\bigr) \cdots f_l\bigl(Z_{t_l}^{\varepsilon}-Z_{t_{l-1}}^\varepsilon
\bigr)\bigr]
\end{eqnarray}
for arbitrary continuous compactly supported functions $f_1,\ldots,f_l\dvtx \R\to\R$.

For fixed set $\IT$, the family of sets $A\in\sigma(Y_t\dvtx t\in\IT)$
for which~(\ref{eq1401-1}) is valid is a Dynkin system provided that
the statement is true for $A=\Omega$. Consequently, it suffices to
prove~(\ref{eq1401-1}) on the $\cap$-stable generator
\[
\cE=\{A_1\cap\cdots\cap A_{l}\dvtx  A_0\in
\cA_0, \ldots, A_{l}\in\cA _l\}, %
\]
where $\cA_1= \sigma(Y_{t_1}),\ldots, \cA_l=\sigma(Y_{t_l}-Y_{t_{l-1}})$.
We note that for $A=A_1\cap\cdots\cap A_l\in\cE$ the random variables
\[
\1_{A_1} f_1\bigl(Z_{t_1}^{n,\varepsilon}\bigr),
\ldots, \1_{A_l} f_l\bigl(Z_{t_l}^{n,\varepsilon}-Z_{t_{l-1}}^{n,\varepsilon}
\bigr)
\]
are independent if
$\IT\subset\varepsilon_n \N_0$ which is fulfilled for sufficiently
large~$n$ since $\IT$ is finite and a subset of $\cD$. Likewise this
holds for $(Z_{t}^{n,\varepsilon})$ replaced by $(Z_t^\varepsilon)$.
Consequently, it suffices to prove that for $k=1,\ldots,l$
\[
\E\bigl[\1_{A_k} f_k\bigl(Z_{t_k}^{n,\varepsilon}-Z_{t_{k-1}}^{n,\varepsilon
}
\bigr)\bigr]\to\E\bigl[\1_{A_k} f_k\bigl(Z_{t_k}^{\varepsilon}-Z_{t_{k-1}}^\varepsilon
\bigr)\bigr]. %
\]
Due to the time homogeneity of the problem, we can and will restrict\break
attention to the case $k=1$ and set $t=t_1$. %
Note\vspace*{1pt} that $\sig(W) \cap\break \bigcup_{\varepsilon'>0} \sigma(\sum_{s\in
(0,t]\dvtx \llvert  \Delta Y_s\rrvert  \geq\varepsilon'} \delta_{\Delta Y_s} )$ is $\cap
$-stable, contains~$\Omega$ and generates a $\sig$-field that
contains $\sigma(Y_t)$.

We conclude that the statement of the lemma is true, if for all $t\in
\cD$, $\varepsilon'>0$, all $A\in\sigma(W)$ and $A'\in\sigma(\sum_{s\in(0,t]\dvtx \llvert  \Delta Y_s\rrvert  \geq\varepsilon'} \delta_{\Delta Y_s})$ and
all continuous compactly supported $f\dvtx \R\to\R$, one has
%
\begin{equation}
\label{eq1411-1} \lim_{n\to\infty} \E\bigl[\1_{A\cap A'} f
\bigl(Z^{n,\varepsilon}_t\bigr)\bigr] = \E\bigl[\1 _{A\cap A'} f
\bigl(Z^\varepsilon_t\bigr)\bigr].
\end{equation}

2\textit{nd step}. In this step, we prove that for $A\in\sig(W)$ and
$A' \in\sig(\Pi)$
\[
\lim_{n\to\infty} \bigl\llvert \E\bigl[\1_{A\cap A'} f
\bigl(Z_t^{n,\varepsilon}\bigr)\bigr]-\P (A) \E\bigl[\1_{A'}
f\bigl(\overline Z_t^{n,\varepsilon}\bigr)\bigr]\bigr\rrvert =0,
\]
where $(\overline Y^{\varepsilon}_s)$ and $(\overline Z_s^{n,\varepsilon})$ are
given by
\[
\overline Y^{\varepsilon}_s=\sigma W_s+ \int
_{(0,s]\times B(0,\varepsilon
)^c} x \,\dd\Pi(u,x) %
\]
and
\[
\overline Z_s^{n,\varepsilon}=\varepsilon_n^{-1/2}
\int_0^s (W_{\iota
_{n+1}(u-)}-W_{\iota_n(u-)})
\,\dd\overline Y^{\varepsilon}_u. %
\]
It suffices to consider the case $\P(A)>0$. We use results of
enlargements of filtrations; see~\cite{jeulinyor}, Theorem~2, page~47,
or~\cite{AnkirchnerDereich}, Example~2: there exists a previsible
process $(\alpha_s)_{s\in[0,T]}$ being\vspace*{1pt} square integrable with respect
to $\P\otimes\ell_{[0,T]}$ such that given $A$ the process
$(W^A_s)_{s\in[0,T]}$
\[
W^A_s:=W_s-\int_0^s
\alpha_u \,\dd u %
\]
is a Wiener process. By Lemma~\ref{le3101-2}, the processes
$(Z_t^{n,\varepsilon})$ and
\[
\overline Z_s ^{n,\varepsilon,A}=\varepsilon_n^{-1/2}
\int_0^s \bigl(W^A_{\iota_{n+1}(u-)}-W^A_{\iota_n(u-)}
\bigr) \,\dd\overline Y^{\varepsilon,A}_u %
\]
with $\overline Y^{\varepsilon,A}=(\sigma W^A_s+\int_{(0,s]\times
B(0,\varepsilon)^c } x \,\dd\Pi(u,x))_{s\in[0,T]}$ are equivalent in
probability. Hence,
\[
\bigl\llvert \E\bigl[\1_{A\cap A'} f\bigl(Z_t^{n,\varepsilon}
\bigr)\bigr]- \E\bigl[\1_{A \cap A'} f\bigl(\overline Z_t^{n,\varepsilon,A}
\bigr)\bigr]\bigr\rrvert \to0.
\]
The set $A$ is independent of $\Pi$. Further, conditionally on $A$ the
process $W^A$ is a Brownian motion that is independent of $\Pi$ which
implies that
\[
\E\bigl[\1_{A \cap A'} f\bigl(\overline Z_t^{n,\varepsilon,A}\bigr)
\bigr]=\P(A) \E\bigl[\1 _{A'} f\bigl(\overline Z_t^{n,\varepsilon}
\bigr)\bigr]. %
\]

3\textit{rd step}.
Let $\Gamma$ denote the finite Poisson point process on
$B(0,\varepsilon')^c$ with
\[
\Gamma=\mathop{\sum_{s\in(0,t]}}_{\llvert  \Delta Y_s\rrvert  \geq\varepsilon'}
\delta_{\Delta Y_s}=\int_{(0,t]\times B(0,\varepsilon)^c} \delta _x \,\dd
\Pi(u,x). %
\]

In the third step, we prove that for every $A'\in\sigma(\Gamma)$ and
every continuous and bounded function $f\dvtx \R\to\R$ one has
\[
\lim_{n\to\infty} \E\bigl[\1_{A'} f\bigl(\overline Z_t^{n,\varepsilon}\bigr)\bigr] =\E\bigl[\1 _{A'} f
\bigl(Z_t^\varepsilon\bigr)\bigr]. %
\]
By dominated convergence, it suffices to show that, almost surely,
%
\begin{equation}
\label{eq0702-1} \lim_{n\to\infty} \E\bigl[ f\bigl(\overline Z_t^{n,\varepsilon}\bigr)\mid \Gamma\bigr] =\E \bigl[f
\bigl(Z_t^\varepsilon\bigr)\mid \Gamma\bigr].
\end{equation}

The regular conditional probability of $\Pi\mid _{(0,t]\times
B(0,\varepsilon')^c}$ given $\Gamma$ can be made precise: the
distribution of $\Pi\mid _{(0,t]\times B(0,\varepsilon')^c}$ given $\{
\Gamma=\gamma:=\sum_{k=1}^m \delta_{y_m}\}$ with $m\in\N$ and
$y_1,\ldots,y_m\in B(0,\varepsilon')^c$ is the same as the
distribution of
\[
\sum_{k=1}^m \delta_{S_k,y_k}
\]
with independent on $(0,t]$ uniformly distributed random variables
$S_1,\ldots,S_m$.
Since, furthermore, $\Pi\mid_{(0,t]\times B(0,\varepsilon')^c}$ is
independent\vspace*{1pt} of $\Pi\mid _{(0,t]\times B(0,\varepsilon')\setminus\{0\}}$
and the Brownian motion $W$, we conclude that the distribution of $\overline Z_t^{n,\varepsilon}$ conditioned on $\{\Gamma=\gamma\}$ equals the
distribution of the random variable
\[
\overline Z_t^{n,\varepsilon,\gamma}= \varepsilon_n^{-1/2}
\int_0^t (W_{\iota^\gamma_{n+1}(u-)}-W_{\iota^\gamma_n(u-)})
\,\dd\overline Y^{n,\gamma}_u %
\]
with $\overline Y^{n,\gamma}_s =\sigma W_s+ \sum_{k=1}^m y_m \1_{\{
\llvert  y_m\rrvert  \geq\varepsilon\}}\1_{\{S_k\leq s\}}$ and
\[
\iota^\gamma_n(s)=\sup \bigl[ \bigl(\varepsilon_n
\Z\cap[0,t]\bigr)\cup\bigl\{ s\in(0,t]\dvtx  h_n\leq\llvert \Delta
Y_s\rrvert <\varepsilon'\bigr\}\cup\{S_1,\ldots,S_m\}\bigr]. %
\]
Here, the random variables $S_1,\ldots,S_m$ are independent of $\Pi
\mid _{(0,t]\times B(0,\varepsilon')}$ and $W$.
Likewise the random variable $Z_t^\varepsilon$ given $\{\Gamma=\gamma
\}$ has the same distribution as the unconditional random variable
\[
Z_t^{\varepsilon,\gamma}= \Upsilon B_t +\sum
_{j=1}^m\frac{\sigma_j}{\sigma} \mathcal
\xi_{j} y_j \1_{\{
\llvert  y_j\rrvert  \geq\varepsilon\}} %
\]
with $\sigma_1,\ldots,\sigma_m$ and $\xi_1,\ldots,\xi_m$ being
independent (also of $B$) with the same distribution as the marks of
the point process $\Pi$.
Consequently, statement\break (\ref{eq0702-1}) follows if for every~$\gamma
$ as above,
\[
\lim_{n\to\infty}\E\bigl[f\bigl(\overline Z_t^{n,\varepsilon,\gamma}
\bigr)\bigr] = \E \bigl[f\bigl(Z_t^{\varepsilon,\gamma}\bigr)\bigr].
\]

We keep $\gamma$ fixed and analyse $\overline Z_t^{n,\varepsilon,\gamma}$
for $n\in\N$ sufficiently large, that is, with $t\in\varepsilon_n\Z$.
We partition $(0,t]$ into $t/\varepsilon_n$ $n$-windows. We call the
$k$th $n$-window to be \emph{occupied} by $S_j$ if~$S_j$ is the only
time in the window $((k-1)\varepsilon_n,k\varepsilon_n]$. Further, we
call a window to be \emph{empty}, if none of the times $S_1,\ldots,S_m$ is in the window.
For each window $k=1,\ldots,t/\varepsilon_n$ that is empty, we set
\[
\cZ_k^{n,\gamma}= \varepsilon_n^{-1/2}
\sigma\int_{(k-1)\varepsilon_n}^{k\varepsilon_n} (W_{\iota^\gamma
_{n+1}(u-)}-W_{\iota^\gamma_n(u-)})
\,\dd W_u, %
\]
and for a window $((k-1)\varepsilon_n,k\varepsilon_n]$ being occupied
by $j$
\[
\cZ_k^{n,\gamma}= \varepsilon_n^{-1/2}
(W_{\iota^\gamma_{n+1}(S_j-)}-W_{\iota^\gamma_n(S_j-)}) y_j \1_{\{
\llvert  y_j\rrvert  \geq\varepsilon\}}.
\]
The remaining $\cZ_k^{n,\gamma}$ can be defined arbitrarily since we
will make use of the fact that the event $\mathcal T_n$ that all
windows are either empty or occupied satisfies $\P(\mathcal T_n)\to
1$.

We first analyse the contribution of the occupied windows. Given that
$\mathcal T_n$ occurs and that $S_1,\ldots,S_m$ are in windows
$k_1,\ldots,k_m$, the random variables $\cZ^{n,\gamma}_{k_1},\ldots,\cZ^{n,\gamma}_{k_m}$ are\vspace*{1pt} independent. We consider their conditional
distributions: conditionally, each $S_j$ is uniformly distributed on
the respective window and the last displacement in $B(0,\varepsilon
')\setminus B(0,h_n)$, respectively, $B(0,h_n)\setminus B(0,h_{n+1})$
has occurred an independent exponentially distributed amount of time ago;
with parameter $\lambda_n=\nu(B(0,\varepsilon')\setminus
B(0,h_n))$, respectively, $\lambda_{n+1}-\lambda_n$. Therefore,
the conditional distribution of $(S_j-\iota_n(S_j), S_j-\iota
_{n+1}(S_j))$ is the same as the one of
\[
\Biggl(\min\bigl(\cU^{\varepsilon_n},\cE^{\lambda_n}\bigr), \sum
_{i=1}^{M} \1 _{((i-1) \varepsilon_n, i\varepsilon_{n}]} \bigl(
\cU^{\varepsilon
_n}\bigr) \min \biggl( \cU^{\varepsilon_n} -\frac{i-1}{M},
\cE^{\lambda
_n}, \cE^{\lambda_{n+1}-\lambda_n} \biggr) \Biggr),
\]
where $\cU^{\varepsilon_n}, \cE^{\lambda_n}$ and $\cE^{\lambda
_{n+1}}$ are independent random variables with $\cU^{\varepsilon_n}$
being uniformly distributed on $[0,\varepsilon_n]$ and $\cE^{\lambda
_n},\cE^{\lambda_{n+1}-\lambda_n}$ being exponentially distributed
with parameters $\lambda_n$ and $\lambda_{n+1}-\lambda_n$.
Consequently, conditionally, one has that
\begin{eqnarray*}
\cZ_{k_j}^{n,\gamma} &\stackrel d=& \varepsilon_n^{-1/2}
\Biggl(\min \bigl(\cU^{\varepsilon_n},\cE^{\lambda_n}\bigr)
\\
&&\hspace*{29pt} {} - \min \Biggl(\sum
_{i=1}^{M} \1 _{((i-1) \varepsilon_n, i\varepsilon_{n}]} \bigl(
\cU^{\varepsilon_n}\bigr) \biggl(\cU^{\varepsilon_n} -\frac{i-1}{M}\biggr),
\cE^{\lambda_n}, \cE ^{\lambda_{n+1}-\lambda_n} \Biggr) \Biggr)^{1/2}
\\
&&{}\times \xi y_j \1_{\{\llvert  y_j\rrvert  \geq\varepsilon\}},
\end{eqnarray*}
where $\xi$ denotes an independent standard normal. By assumption,
$\lambda_n/\varepsilon_n\to\theta$ as $n\to\infty$ so that the
latter distribution converges to the one of $\frac{\sigma_j}{\sigma}
\xi_j y_j$. Hence, conditionally on $\mathcal T_n$ one has
\[
\mathop{\sum_{k\in\N\cap[0,t/\varepsilon_n]}}_{k\mathrm{th}~n\mbox{-}\mathrm{window}~\mathrm{occupied}}
\cZ_k^{n,\gamma}\Rightarrow\sum_{j=1}^m
\frac
{\sigma_j}{\sigma} \xi_j y_j \1_{\{\llvert  y_j\rrvert  \geq\varepsilon\}}.
\]

Next, we analyse the contribution of all empty windows. Given $\mathcal
T_n$, there are $t/\varepsilon_n-m$ empty windows and the
corresponding random variables $\cZ_k^{n,\gamma}$ are independent and
identically distributed. We have
\[
\E\bigl[\cZ^{n,\gamma}_1\mid (0,\varepsilon_n]
\mbox{ empty, }\cT_n\bigr]=0 %
\]
since $W$ is independent of the event we condition on. Further, by It\^
o's isometry and the scaling properties of Brownian motion one has
%
\begin{eqnarray}\label{eq2602-6}
&& \Var \bigl(\cZ^{n,\gamma}_1 \mid (0,
\varepsilon_n]\mbox{ empty, }\cT_n \bigr)\nonumber
\\
&&\qquad =
\varepsilon_n^{-1}\sigma^2 \E \biggl[ \int
_0^{\varepsilon_n} (W_{\iota^\gamma_{n+1}(u)}-W_{\iota^\gamma
_{n}(u)})^2
\,\dd u \mid (0,\varepsilon_n]\mbox{ empty, }\cT_n \biggr]
\nonumber\\[-8pt]\\[-8pt]\nonumber
&&\qquad =\varepsilon_n \sigma^2 \E \bigl[ (W_{\varepsilon_n^{-1}\iota
^\gamma_{n+1}(\cU^{\varepsilon_n})}-W_{\varepsilon_n^{-1}\iota
^\gamma_{n}(\cU^{\varepsilon_n})})^2
\mid (0,\varepsilon _n]\mbox{ empty, }\cT_n \bigr]
\\
&&\qquad = \varepsilon_n \sigma^2 \E \bigl[
\varepsilon_n^{-1}\iota ^\gamma_{n+1}
\bigl(\cU^{\varepsilon_n}\bigr)-\varepsilon_n^{-1}
\iota^\gamma _{n}\bigl(\cU^{\varepsilon_n}\bigr)\bigr] \mid (0,
\varepsilon_n]\mbox{ empty, }\cT_n \bigr].\nonumber
\end{eqnarray}
Here, we denote again by $\cU^{\varepsilon_n}$ an independent uniform
random variable on $[0,\varepsilon_n]$ and we used that conditionally
the processes $\iota^\gamma_n$ and~$\iota^\gamma_{n+1}$ are
independent of the Brownian motion~$W$.
As above, we note that the distributions of $ \varepsilon_n^{-1}\iota
^\gamma_{n+1}(\cU^{\varepsilon_n})$ and $\varepsilon_n^{-1}\iota
^\gamma_{n}(\cU^{\varepsilon_n})$ are identically distributed as
\[
\varepsilon_n^{-1} \biggl(\cE^{\lambda_{n+1}}\wedge
\frac{\cU
^{\varepsilon_n}}M \biggr)\quad\mbox{and}\quad\varepsilon _n^{-1}
\bigl(\cE^{\lambda_{n}}\wedge\cU^{\varepsilon_n} \bigr). %
\]
By assumption \textup{(ML2)}, these converge in $L^1$ to $\cE^{M\theta}\wedge
\cU^{1/M} $ and $\cE^{\theta}\wedge\cU^{1}$, respectively. Hence,
computing the respective expectations gives with~(\ref{eq2602-6})
\[
\varepsilon_n^{-1} \Var \bigl(\cZ^{n,\gamma}_1
\mid (0,\varepsilon_n]\mbox{ empty, }\cT_n \bigr)\to
\sigma^2 \frac
{M-1}{M} \frac{e^{-\theta}-(1-\theta)}{\theta^2}=:\Upsilon^2.
\]
The uniform $L^2$-integrability of $\cL(\varepsilon_n^{-1/2} \cZ
^{n,\gamma}_1 \mid  (0,\varepsilon_n]\mbox{ empty, }\cT_n)$ follows
by noticing that by the Burkh\"older--Davis--Gundy inequality there
exists a universal constant $\kappa$ such that
\begin{eqnarray*}
&& \E \bigl(\bigl(\cZ^{n,\gamma}_1\bigr)^4 \mid (0,
\varepsilon_n]\mbox{ empty, }\cT_n \bigr)
\\
&&\qquad \leq\kappa
\varepsilon_n^{-2}\sigma^4 \E \biggl[ \biggl(
\int_0^{\varepsilon_n} (W_{\iota^\gamma
_{n+1}(u)}-W_{\iota^\gamma_{n}(u)})^2
\,\dd u \biggr)^2 \Big| (0,\varepsilon_n] \mbox{
empty, }\cT_n \biggr]
\\
&&\qquad \leq  4\kappa \sigma^4 \E \Bigl[\sup_{u\in[0,\varepsilon_n]}
W_{u}^4 \Bigr]=4\kappa \sigma^4
\varepsilon_n^{2} \E \Bigl[\sup_{u\in[0,1]}
W_{u}^4 \Bigr].
\end{eqnarray*}
Hence, conditionally on $\mathcal T_n$ one has
\[
\mathop{\sum_{k\in\N\cap(0,t/\varepsilon_n]}}_{k\mathrm{th}~n\mbox{-}\mathrm{window}~\mathrm{empty}}
\cZ_k^{n,\gamma}\Rightarrow\cN\bigl(0,\Upsilon^2 t
\bigr). %
\]

Given $\mathcal T_n$ the contribution of the empty and occupied windows
are independent, so that since $\P(\mathcal T_n)\to1$, generally
\[
\sum_{k=1}^{t/\varepsilon_n} \cZ_k^{n,\gamma}
\Rightarrow Z_t^{\varepsilon,\gamma}. %
\]
It remains to show that
\[
\lim_{n\to\infty} \Biggl(\overline Z^{n,\varepsilon,\gamma}_t- \sum
_{k=1}^{t/\varepsilon_n} \cZ_k^{n,\gamma}
\Biggr)=0\qquad\mbox {in probability.} %
\]
This follows immediately by noticing that, given $\mathcal T_n$, one has
\begin{eqnarray*}
&& \overline Z^{n,\varepsilon,\gamma}- \sum_{k=1}^{t/\varepsilon_n}
\cZ_k^{n,\gamma}
\\
&&\qquad =\sigma\varepsilon _n^{-1/2}
\mathop{\sum_{k\in\N\cap[0,t/\varepsilon_n]}}_{k\mathrm{th}~n\mbox{-}\mathrm{window}~\mathrm{occupied}}\int
_{(k-1)\varepsilon_n}^{k\varepsilon
_n}(W_{\iota_{n+1}(u-)}-W_{\iota_n(u-)})
\,\dd W_u, %
\end{eqnarray*}
where the sum on the right-hand side is over $m$ independent and
identically distributed summands each having second moment smaller than
$\varepsilon_n^2$.

4\textit{th step}. In the last step, we combine the results of the
previous steps.
By step one, it suffices to verify equation~(\ref{eq1411-1}). Provided
that the statement is true for $A=\Omega$, the system of sets $A$ for
which~(\ref{eq1411-1}) is satisfied is a Dynkin system. Consequently,
it suffices to verify validity for sets $A\cap A'$ with $A\in\sigma
(W_t)$ and $A'\in\sigma(\Gamma)$. By step two, one has
\[
\lim_{n\to\infty}\bigl\llvert \E[\1_{A\cap A'} f
\bigl(Z^{n,\varepsilon}_t\bigr)- \P (A) \E\bigl[\1_{A'} f
\bigl(\overline Z^{n,\varepsilon}_t\bigr)\bigr]\bigr\rrvert \to0 %
\]
and by step three
\[
\lim_{n\to\infty} \E\bigl[\1_{A'} f\bigl(\overline Z^{n,\varepsilon}\bigr)\bigr]=\E\bigl[\1 _{A'} f\bigl(Z_t^\varepsilon
\bigr)\bigr] %
\]
so that
\[
\lim_{n\to\infty} \E[\1_{A\cap A'} f\bigl(Z^{n,\varepsilon}_t
\bigr)= \P(A) \E\bigl[\1_{A'} f\bigl(Z_t^\varepsilon
\bigr)\bigr]. %
\]
The proof is complete by noticing that $\sig(W_t)$ is independent of
$\sigma(\Gamma, Z_t^\varepsilon)$ so that
\[
\P(A) \E\bigl[\1_{A'}f\bigl(\overline Z^\varepsilon_t
\bigr)\bigr]=\E\bigl[\1_{A\cap A'}f\bigl(\overline Z^\varepsilon_t
\bigr)\bigr].
\]\upqed
\end{pf}

\section{Scaled errors of derived quantities}\label{derivedquantities}

In this section, we collect results that will enable us to deduce the
main central limit theorems with the help of Theorem~\ref{twosuccessive}.

\subsection{The integrated processes}
The following lemma is central to the proof of Theorems~\ref
{theo2402-2} and~\ref{theo2003-1}.

\begin{lemma}\label{le2103-2}
If assumptions \textup{(ML1)} and \textup{(ML2)} hold, then one has
\[
\lim_{n\to\infty} \varepsilon_n^{-1} \E
\biggl[ \biggl\llvert \int_0^T \bigl(\widehat{X}^n_t-\overline X^n_t\bigr) \,\dd t
\biggr\rrvert ^2 \biggr]=0. %
\]
\end{lemma}

\begin{pf}
With $b_n:=b-\int_{B(0,h_n)^c}x \nu(\dd x)$ we have for
$t\in[0,T]$
%
\begin{eqnarray}\label{eq1702-1}
\widehat{X}^n_t-\overline X^n_t&=&
a(\widehat{X}_{\iota
_n(t)}) \bigl(Y_{t}^h-Y_{\iota_n(t)}^h
\bigr)
\nonumber\\[-8pt]\\[-8pt]\nonumber
&=& a(\widehat{X}_{\iota_n(t)}) \bigl(b_n \bigl(t-
\iota_n(t)\bigr)+\sigma(W_{t}-W_{\iota_n(t)})\bigr).
\end{eqnarray}
We estimate
\[
\E \biggl[\biggl\llvert \int_0^T a(\widehat{X}_{\iota_n(t)}) b_n \bigl(t-\iota _n(t)\bigr) \,\dd t
\biggr\rrvert ^2 \biggr] \leq b_n^2
\varepsilon_n^2 T \E \biggl[\int_0^T
\bigl\llvert a\bigl(\widehat{X}^n_{\iota_n(t-)}\bigr)\bigr\rrvert
^2 \,\dd t \biggr].
\]
The latter expectation is uniformly bounded over all $n$; see
Lemma~\ref{boundapp}. Further, $b_n^2=o(\varepsilon_n^{-1})$ by
Lemma~\ref{driftdisappear}. Consequently, the first term is of order
$o(\varepsilon_n)$.
By Fubini,
\begin{eqnarray*}
&& \E \biggl[ \biggl(\int_0^T a(\widehat{X}_{\iota_n(t)}) \sigma (W_t-W_{\iota_n(t)}) \,\dd t
\biggr)^2 \biggr]
\\
&&\qquad =\sigma^2 \int_0^T\int
_0^T \E\bigl[a(\widehat{X}_{\iota_n(t)})
(W_t-W_{\iota_n(t)}) a(\widehat{X}_{\iota_n(u)})
(W_u-W_{\iota_n(u)}) \bigr] \,\dd t \,\dd u.
\end{eqnarray*}
Further, for $0\leq t\leq u\leq T$,
\begin{eqnarray*}
&& \E\bigl[a(\widehat{X}_{\iota_n(t)}) (W_t-W_{\iota_n(t)}) a(\widehat{X}_{\iota
_n(u)}) (W_u-W_{\iota_n(u)})\mid
\iota_n, \widehat{X}_{\iota_n(t)}\bigr]
\\
&&\qquad =\1_{\{
\iota_n(t)=\iota_n(u)\}} a(\widehat{X}_{\iota_n(t)})^2 \bigl((t\wedge u)-\iota_n(t)\bigr)
\end{eqnarray*}
and since the statement is symmetric in the variables $t,u$ also for
$0\leq u\leq t\leq T$. Consequently,
\[
\E \biggl[ \biggl(\int_0^T a(\widehat{X}_{\iota_n(t)}) \sigma (W_t-W_{\iota_n(t)}) \,\dd t
\biggr)^2 \biggr] \leq2\varepsilon_n^2
\sigma^2 \int_0^T \E\bigl[a(\widehat{X}_{\iota_n(t)}) ^2\bigr] \,\dd t.
\]
We recall that the latter expectation is uniformly bounded so that this
term is also of order $o(\varepsilon_n)$.
\end{pf}

\subsection{The supremum}

The results of this subsection are central to the proof of Theorem~\ref{thm2602-1}.
We first give some qualitative results for solutions $X=(X_t)_{t\in
[0,T]}$ of the stochastic differential equation
\[
\dd X_t= a(X_{t-}) \,\dd Y_t %
\]
with arbitrary starting value. We additionally assume that $a$ does
\emph{not} attain zero.

\begin{lemma}\label{le1302-2}
One has for every $t\in[0,T]$ that, almost surely,
\[
\sup_{s\in[0,t]} X_s> X_0\vee
X_t. %
\]
%
\end{lemma}

\begin{pf}We only prove that
\[
\sup_{s\in[0,t]} X_s> X_t %
\]
and remark that the remaining statement follows by similar simpler
considerations.

\begin{longlist}[1\textit{st step}.]
\item[1\textit{st step}.] In the first step, we show that
\[
\frac{1}{{\sqrt\varepsilon}}(X_{t-\varepsilon+\varepsilon
s}-X_{t-\varepsilon})_{s\in[0,1]}
\stackrel{\operatorname {stably}} {\Longrightarrow}\bigl(\sigma a(X_t)
B_s\bigr)_{s\in[0,1]}. %
\]
We show the statement in two steps: first note that
\[
\frac{1}{{\sqrt\varepsilon}}(X_{t-\varepsilon+\varepsilon
s}-X_{t-\varepsilon})_{s\in[0,1]}\quad
\mbox{and}\quad\frac{1}{{\sqrt\varepsilon}} \bigl(a(X_{t-\varepsilon}) (Y_{t-\varepsilon
+\varepsilon s}-Y_{t-\varepsilon})
\bigr)_{s\in[0,1]} %
\]
are equivalent in ucp. Further, $Z^\varepsilon:=(\varepsilon^{-1/2}
(Y_{T-\varepsilon+\varepsilon s}-Y_{t-\varepsilon}))_{s\in[0,1]}$ is
independent of $a(X_{t-\varepsilon})$ and $a(X_{t-\varepsilon})$
tends to $a(X_t)$, almost surely.
Hence, it remains to show that $Z^\varepsilon$ converges for
$\varepsilon\downarrow0$ in distribution to $\sigma B$.
Note that $Z^\varepsilon$ is a L\'evy-process with triplet $(b\sqrt
\varepsilon,\sigma^2,\nu_\varepsilon)$, where $\nu_\varepsilon
(A)=\varepsilon \nu(\sqrt\varepsilon A)$ for Borel sets $A\subset
\R\setminus\{0\}$.
It suffices to show that L\'evy-processes $\overline Z^\varepsilon$ with
triplet $(0,0,\nu_\varepsilon)$ converge to the zero process.

We uniquely represent $\overline Z^\varepsilon$ as
\[
\overline Z^\varepsilon_t = \overline Z^{\varepsilon,r}_t+
\hspace*{1pt}\overline{\hspace*{-1pt}\overline Z}^{\hspace*{1pt}\varepsilon,r}- b_{\varepsilon,r} t %
\]
with independent L\'evy processes $\overline Z^{\varepsilon,r}_t$ and $\hspace*{1pt}\overline{\hspace*{-1pt}\overline Z}^{\hspace*{1pt}\varepsilon,r}$, the first one with triplet $(0,0,\nu
_\varepsilon\mid _{B(0,r)})$, the second one being a compound Poisson
process with intensity $\nu\mid _{B(0,r)^c}$, and with $b_{\varepsilon,r}:=\int_{B(0,r)^c} x  \,\dd\nu_\varepsilon(x)$. Clearly, for
$\delta>0$
%
\begin{equation}\label{eq1402-1}
\P \Bigl(\sup_{t\in[0,1]} \bigl\llvert \overline Z^\varepsilon_t\bigr\rrvert >\delta \Bigr)\leq
\1_{\{\llvert  b_{\varepsilon,r}\rrvert  >\delta/2\}}+ \P \Bigl(\sup_{t\in[0,1]} \bigl\llvert \overline Z^{\varepsilon,r}_t\bigr\rrvert >\delta/2 \Bigr)+ \P\bigl(
\hspace*{1pt}\overline{\hspace*{-1pt}\overline Z}^{\hspace*{1pt}\varepsilon,r}\neq0\bigr).\hspace*{-30pt}
\end{equation}
For $r>0$, one has
\begin{eqnarray*}
r \nu\bigl(B(0,r)^c\bigr)&\leq&\int_{B(0,r)^c} \llvert
x\rrvert \,\dd\nu_\varepsilon (x)= \varepsilon\int_{B(0,\sqrt\varepsilon r)^c}
\frac{\llvert  x\rrvert  }{\sqrt
\varepsilon} \,\dd\nu(x)
\\
& \leq& \sqrt\varepsilon \int_{B(0,\sqrt\varepsilon r)} \frac
{x^2}{\sqrt\varepsilon r} \nu(\dd x) \leq\frac{1}{r} \int x^2 \nu (\dd x).
\end{eqnarray*}
Hence, $\llvert  b_{\varepsilon,r}\rrvert  \leq\delta/2$, for sufficiently large
$r$, and $\P(\overline{\overline Z}^{\varepsilon,r}\neq0)\leq\nu
(B(0,r)^c)\leq\frac{1}{r^2} \int x^2 \nu(\dd x)$. Further,
\[
\int_{B(0,r)} x^2 \,\dd\nu_\varepsilon(x)= \int
_{B(0,\sqrt
\varepsilon r)} x^2 \,\dd\nu(x) \to0 %
\]
so that Doob's $L^2$-inequality yields
\[
\lim_{\varepsilon\downarrow0} \P \Bigl(\sup_{t\in[0,1]} \bigl
\llvert \overline Z^{\varepsilon,r}_t\bigr\rrvert >\delta/2 \Bigr)=0.
\]
Plugging these estimates into~(\ref{eq1402-1}) gives
\[
\limsup_{\varepsilon\downarrow0} \P \Bigl(\sup_{t\in[0,1]} \bigl
\llvert \overline Z^\varepsilon_t\bigr\rrvert >\delta \Bigr)\leq
\frac{1}{r^2} \int x^2 \nu(\dd x) %
\]
and the statement of step one follows by noticing that $r>0$ can be
chosen arbitrarily large.

\item[2\textit{nd step}.] Clearly, for $\varepsilon\in(0,t]$,
\[
\P\Bigl(\sup_{s\in[0,t]} X_s= X_t\Bigr)
\leq\P \Bigl(\varepsilon^{-1/2}\sup_{s\in[0,1]}(X_{t-\varepsilon+\varepsilon s}-X_{t-\varepsilon})=
\varepsilon^{-1/2} (X_t-X_{t-\varepsilon}) \Bigr).
\]
The set of all c\`adl\`ag functions $x\dvtx [0,1]\to\R$ with $\sup_{s\in
[0,1]} x_s= x_1$
is closed in the Skorokhod space so that
\begin{eqnarray*}
&& \P \Bigl(\sup_{s\in[0,t]} X_s= X_t \Bigr)
\\
&&\qquad \leq\limsup_{\varepsilon
\downarrow0} \P \Bigl( \varepsilon^{-1/2}\sup
_{s\in
[0,1]}(X_{t-\varepsilon+\varepsilon s}-X_{t-\varepsilon})=
\varepsilon^{-1/2} (X_t-X_{t-\varepsilon}) \Bigr)
\\
&&\qquad \leq \P \Bigl( a(X_t) \sup_{s\in[0,1]} \sigma
B_s= a(X_t) B_1 \Bigr)=0.
\end{eqnarray*}\vspace*{-28pt}
\end{longlist}
\end{pf}

\begin{lemma}\label{le2502-3} Suppose that $a(x)\neq0$ for all $x\in
\R$. There is a unique random time $S$ (up to indistinguishability)
such that, almost surely,
\[
\sup_{s\in[0,T]} X_s= X_S %
\]
and one has $\Delta X_S=0$.
Further, for every $\varepsilon>0$, almost surely,
\[
\sup_{s\in[0,S]\dvtx  \llvert  s-S\rrvert  \geq\varepsilon} X_s< X_S. %
\]
\end{lemma}

\begin{pf}
1\textit{st step}. First we prove that the supremum $\sup_{t\in[0,T]}X_t$ is almost surely attained at some random time $S$
with $\Delta X_S=0$. By compactness of the time domain, we can find an
almost surely convergent $[0,T]$-valued sequence $(S_n)_{n\in\N}$ of
random variables, say with limit $S$, with
\[
\lim_{n\to\infty} X_{S_n}= \sup_{t\in[0,T]}
X_t. %
\]
%

Let $h>0$. We represent $Y$ as sum
\[
Y_t=Y^h_t+\sum_{k=1}^N
\1_{[T_k,T]}(t) \Delta Y_{T_k}, %
\]
where $T_1,\ldots,T_N$ are the increasingly ordered times of the
discontinuities of $Y$ being larger than $h$. Further, $Y^h$ is a L\'
evy process that is independent of $\overline Y^h:=Y-Y^h$. Given $\overline Y^h$,
for every $k=1,\ldots,N$, the process $(X_t)_{t\in[T_{k-1},T_k)}$
solves the~SDE
\[
\dd X_t=a(X_{t-}) \,\dd Y^h_t
\]
and we have, almost surely, that
\[
X_{T_k-}= X_{T_{k-1}} +\int_{T_{k-1}}^{T_k}
a(X_s) \,\dd Y^h_s. %
\]
Consequently, we can apply Lemma~\ref{le1302-2} and conclude that,
almost surely, for each $k=1,\ldots,N+1$,
\[
\sup_{s\in[T_{k-1},T_{k})} X_s > X_{T_{k-1}} \vee
X_{T_k-} %
\]
with $T_0=0$ and $T_{N+1}=T$. Hence, almost surely,
\[
\sup_{s\in[0,T] } X_s> \sup_{k=1,\ldots,N+1}
X_{T_{k-1}} \vee X_{T_k-}. %
\]
Consequently, $S$ is almost surely not equal to $0$ or $T$ or a time
with displacement larger than $h$. Since $h>0$ was arbitrary, we get
that, almost surely, $\Delta X_S=0$, so that
\[
X_S=\lim_{n\to\infty} X_{S_n}=\sup
_{t\in[0,T]} X_t\qquad\mbox {almost surely.} %
\]

2\textit{nd step}. We prove that for every $t\in[0,T]$ the distribution
of $\sup_{s\in[0,t]} X_s$ has no atom. Suppose that it has an atom in
$z\in\R$. We consider the stopping time
\[
T_{\{z\}}=\inf\bigl\{t\in[0,T]\dvtx  X_t =z\bigr\} %
\]
with the convention $T_{\{z\}}=\infty$ in the case when $z$ is not
hit. For $\varepsilon>0$, conditionally on the event $\{T_{\{z\}}\leq
T-\varepsilon\}$ the process $(\widetilde X_{s})_{s\in[0,\varepsilon]}$ with
\[
\widetilde X_s= X_{T_{\{z\}}+s} %
\]
starts in $z$ and solves $\dd\widetilde X_s= a(\widetilde X_s) \,\dd\widetilde
Y_s$ with $\widetilde Y$ denoting the $T_{\{z\}}$-shifted L\'evy
process~$Y$. Hence, by Lemma~\ref{le1302-2}, one has almost surely on
$\{T_{\{z\}}\leq T-\varepsilon\}$ that
\[
z= \widetilde X_0< \sup_{s\in[0,\varepsilon]} \widetilde
X_s\leq\sup_{s\in
[0,T]} X_s. %
\]
Since $\varepsilon>0$ is arbitrary and $X$ does not attain its
supremum in $T$, it follows that $\P( \sup_{s\in[0,T]}
X_s=z)=0$.

3\textit{rd step}. We prove that the supremum over two disjoint time
windows $[u,v)$ and $[w,z)$ with $0\leq u<v\leq w<z\leq T$, satisfies
\[
\sup_{s\in[u,v)} X_s \neq\sup_{s\in[w,z)}
X_s, %
\]
almost surely. By the Markov property, the random variables $\sup_{s\in[u,v)} X_s$ and $\sup_{s\in[w,z)} X_s$ are independent given
$X_w$ and we get
\begin{eqnarray*}
&& \P \Bigl(\sup_{s\in[u,v)} X_s
= \sup
_{s\in[w,z)} X_s \Bigr)
\\
&&\qquad= \int\P \Bigl(\sup
_{s\in[w,z)} X_s=y\big| X_w=x \Bigr) \,\dd\P
_{(X_w,\sup_{s\in[u,v)} X_s )}(x,y), %
\end{eqnarray*}
were $ \P_{(X_w,\sup_{s\in[u,v)} X_s )}$ denotes the distribution of
$(X_w,\sup_{s\in[u,v)} X_s )$. We note that the conditional process
$(X_s)_{s\in[w,z)}$ is again a solution of the SDE started in $x$ and
by step two the inner conditional probability equals zero.

4\textit{th step}. We finish the proof of the statement. For given
$\varepsilon>0$, we choose deterministic times $0=t_0<t_1<\cdots
<t_m=T$ with $t_k-t_{k-1}\leq\varepsilon$. By step three, there is,
almost surely, one window in which the supremum is attained, say in
$[t_{M-1},t_M)$, and
\[
\sup_{s\in[0,T]\dvtx  \llvert  S-s\rrvert  \geq\varepsilon} X_s\leq\sup_{k\in\{
1,\ldots,m\}\setminus\{M\}}
\sup_{s\in[t_{k-1},t_k)} X_s< \sup_{s\in[t_{M-1},t_M)}
X_s=X_S. %
\]\upqed
\end{pf}

\begin{lemma}\label{le2502-4}Suppose that $a(x)\neq0$ for all $x\in
\R$ and denote by $S$ the random time at which $X$ attains its
maximum. One has
\[
\varepsilon_n^{-1/2} \Bigl( \sup_{t\in[0,T]}
X_t^{n+1}- \sup_{t\in
[0,T]}
X^{n}_t \Bigr)-U^{n,n+1}_S\to0\qquad
\mbox{in probability}. %
\]
\end{lemma}

\begin{pf}
With Lemma~\ref{le2502-3} we conclude that, for every $\varepsilon
>0$, one has with high probability that
\begin{eqnarray*}
&& \Bigl\llvert \varepsilon_n^{-1/2} \Bigl(\sup
_{t\in[0,T]} X_t^{n+1}- \sup
_{t\in[0,T]} X^{n}_t \Bigr)-U^{n,n+1}_S
\Bigr\rrvert
\\
&&\qquad \leq\sup_{t\dvtx  \llvert  t-S\rrvert
\leq\varepsilon} \bigl\llvert \varepsilon_n^{-1/2}
\bigl( X_t^{n+1}- X_t^n
\bigr)-U_S^{n,n+1}\bigr\rrvert
\\
&&\qquad = \sup_{t\dvtx \llvert  t-S\rrvert  \leq\varepsilon} \bigl\llvert U^{n,n+1}_t-U_S^{n,n+1}
\bigr\rrvert.
\end{eqnarray*}
For $\varepsilon,\delta>0$, consider
\[
A_{\varepsilon,\delta}=\Bigl\{(s,x)\in[0,T]\times\ID(\R)\dvtx  \sup_{(t,u)\dvtx s-\varepsilon\leq t\leq u\leq s+\varepsilon}\dvtx
\llvert x_t-x_u\rrvert \geq \delta\Bigr\}. %
\]
Note that $\mathrm{cl}(A_{\varepsilon,\delta})\subset
A_{2\varepsilon,\delta}$ and recall that $(S,U^{n,n+1})\Rightarrow
(S,U)$. Hence,
\begin{eqnarray*}
&& \limsup_{n\to\infty}  \P \Bigl(\Bigl\llvert \varepsilon_n^{-1/2}
\sup_{t\in[0,T]} X_t^{n+1}-
\varepsilon_n^{-1/2} \sup_{t\in[0,T]}
X^{n}_t-U^{n,n+1}_S\Bigr\rrvert \geq
\delta \Bigr)
\\
&&\qquad \leq\limsup_{n\to\infty} \P\bigl(\bigl(S,U^{n,n+1}\bigr)\in
A_{\varepsilon,\delta}\bigr)\leq\P\bigl((S,U)\in A_{2\varepsilon,\delta}\bigr).
\end{eqnarray*}
Note that $U$ is almost surely continuous in $S$ so that for
$\varepsilon\downarrow0$, $\P((S,U)\in  A_{2\varepsilon,\delta})\to0$.
\end{pf}

\section{Proofs of the central limit theorems}\label{proofCLT}

In this section, we prove all central limit theorems and Theorem~\ref
{theo0303-1}. We will verify the Lindeberg conditions for the summands
of the multilevel estimate $\widehat S(F)$; see~(\ref
{multiestimator}). As shown in Lemma~\ref{le1702-4} in the \hyperref[append]{Appendix}, a
central limit theorem holds for the idealised approximations
$X^1,X^2,\ldots,$ if:
\begin{longlist}[(2)]
\item[(1)] $\lim_{n\to\infty} \Var(\varepsilon_{n}^{-1/2}
(F(X^{n+1})-F(X^{n}))= \rho^2$ and\vspace*{1pt}

\item[(2)] $(\varepsilon_{n}^{-1/2} (F(X^{n+1})-F(X^{n}))\dvtx k\in\N)$
is uniformly $L^2$-integrable.
\end{longlist}
The section is organised as follows. In Section~\ref{sec41}, we
verify uniform $L^2$-in\-tegrability of the error process in supremum
norm which will allow us to verify property~(2) in the central limit
theorems. In Section~\ref{proofCLT1}, we prove Theorems~\ref
{theo2402-1} and~\ref{theo2402-2}, essentially by verifying property~(1).

It remains to deduce Theorems~\ref{thm2602-1} and~\ref{theo2003-1}
from the respective theorems for the idealised scheme. By Lemmas~\ref
{le2602-1},~\ref{le2602-2} and~\ref{le2103-2}, switching from the
idealised to the continuous or piecewise constant approximation leads
to asymptotically equivalent $L^2$-errors. Hence, the same error
process can be used and, in particular, uniform $L^2$-integrability
prevails due to Lemma~\ref{lemmauniforminteg}. Consequently, the
identical proofs yield the statements.

Finally, we prove Theorem~\ref{theo0303-1} in Section~\ref{sec43}.

\subsection{Uniform $L^2$-integrability}\label{sec41}

\begin{proposition}\label{prop2502-1}
The sequence $(\varepsilon_n^{-1/2} \sup_{t\in
[0,T]}\llvert  X^{n+1}_t-X^n_t\rrvert  )_{n\in\N}$
is uniformly $L^2$-integrable.
\end{proposition}

To prove the proposition, we will make use of the perturbation
estimates given in the \hyperref[append]{Appendix}; see Section~\ref{secperturb}.
Recall that $U^{n,n+1}=\varepsilon_n^{-1/2} (X^{n+1}-X^n)$ satisfies
the equation
\begin{eqnarray*}
U^{n,n+1}_{t}&=&\int_0^t
D_{s-}^{n,n+1} U^{n,n+1}_{s-} \,\dd
Y_{s}+\varepsilon_n^{-1/2} \int
_0^t D^n_{s-}
A_{s-}^n (Y_{s-}-Y_{\iota_{n}(s-)}) \,\dd
Y_{s}
\\
&&{}{} -\varepsilon_n^{-1/2} \int_0^t
D_{s-}^{n+1} A_{s-}^{n+1}
(Y_{s-}-Y_{\iota_{n+1}(s-)}) \,\dd Y_{s}.
\end{eqnarray*}
We use approximations indexed by $m\in\N$: we denote by
\[
\cU
^{n,n+1,m}=\bigl(\cU^{n,n+1,m}_t\bigr)_{t\in[0,T]}
\]
the solution of the equation
%
\begin{eqnarray}
\label{ggl02}
\cU^{n,n+1,m}_t&=&\int_0^tD^{n,n+1}_{s-}
\cU^{n,n+1,m}_s \,\dd\cY ^m_s\nonumber
\\
&&{} + \varepsilon_n^{-1/2}\sigma\int_0^t
D^{n}_{s-}\cA ^{n,m}_{s-}(W_{s-}-W_{\iota_{n}(s-)})
\,\dd\cY^m_s
\\
&&{}-\varepsilon_n^{-1/2}\sigma\int_0^t
D^{n+1}_{s-}\cA ^{n,m}_{s-}(W_{s-}-W_{\iota_{n+1}(s-)})
\,\dd\cY^m_s,\nonumber
\end{eqnarray}
where $\cY^m=(\cY^m_t)_{t\in[0,T]}$ is given by
\[
\cY^m_t=bt+\sigma W_t+\lim
_{\delta\downarrow0} \int_{(0,t]\times
(B(0,m)\setminus B(0,\delta))} x \,\dd\overline\Pi(s,x),
\]
and $\cA^{n,m}=(\cA^{n,m}_t)_{t\in[0,T]}$ is the simple adapted
c\`adl\`ag process given by
\[
\cA^{n,m}_t =\cases{A^n_t, &\quad
if $\bigl\llvert A^n_t\bigr\rrvert \leq m$,
\cr
0, &\quad
else.} %
\]
The proof of the proposition is achieved in two steps. We show that:
\begin{longlist}[2.]
\item[1.] $\lim_{m\uparrow\infty} \limsup_{n\to\infty} \E[\sup_{t\in[0,T]} \llvert  U^{n,n+1}_t-\cU^{n,n+1,m}_t\rrvert  ^2]=0$ and

\item[2.] for every $p\geq2$ and $m\in\N$, $\E[\sup_{t\in[0,T]} \llvert  \cU
^{n,n+1,m}_t\rrvert  ^p]<\infty$.
\end{longlist}
Then the uniform $L^2$-integrability of $(\sup_{t\in[0,T]}
\llvert  U^{n,n+1}_t\rrvert  )_{n\in\N}$ follows with Lemma~\ref
{lemmauniforminteg}.

\begin{lemma}
One has
\[
\lim_{m\uparrow\infty}\limsup_{n\to\infty}\E \Bigl[\sup
_{t\in
[0,T]}\bigl\llvert U^{n,n+1}_t-
\cU^{n,n+1,m}_t\bigr\rrvert ^2 \Bigr]=0. %
\]
\end{lemma}

\begin{pf}The processes $\cU^{n,n+1,m}$ are perturbations of
$U^{n,n+1}$ as analysed in Lemma~\ref{estimatesdes2}. More explicitly,
the result follows if there exists a constant $\kappa>0$ such that
%
\begin{equation}
\label{ggl03} \E \biggl[\sup_{t\in[0,T]}\biggl\llvert
\varepsilon_n^{-1/2}\int_0^t
D^n_{s-}\cA^{n,m}_{s-}(W_{s-}-W_{\iota_n(s-)})
\,\dd\cY^m_s\biggr\rrvert ^2 \biggr]\leq
\kappa,
\end{equation}
for all $n,m\in\N$, and
%
\begin{eqnarray}\label{ggl05}
\qquad&& \lim_{m\to\infty}\limsup_{n\to\infty}
\varepsilon_n^{-1} \E \biggl[\sup_{t\in[0,T]}
\biggl\llvert \int_0^tD^n_{s-}A^n_{s-}(Y_{s-}-Y_{\iota_{n}(s-)})
\,\dd Y_s
\nonumber\\[-8pt]\\[-8pt]\nonumber
&&\hspace*{116pt}{}-\sigma\int_0^tD^n_{s-}
\cA^{n,m}_{s-}(W_{s-}-W_{\iota_{n}(s-)})\,\dd
\cY^m_s\biggr\rrvert ^2 \biggr]=0.
\end{eqnarray}
Using Lemma~\ref{estimateintegral}, the uniform boundedness of $D^n$,
conditional independence of $\cA^{n,m}_{s-}$ and $W_{s-}- W_{\iota
_{n}(s-)}$ given $\iota_n$, there exists a constant $\kappa_1>0$ such that
\begin{eqnarray*}
&& \varepsilon_n^{-1} \E \biggl[\sup_{0\leq r\leq t}
\biggl\llvert \int_0^r D^{n}_{s-}
\cA^{n,m}_{s-}(W_{s-}-W_{\iota_{n}(s-)}) \,\dd
\cY^m_s\biggr\rrvert ^2 \biggr]
\\
&&\qquad \leq
\kappa_1 \int_0^T
\varepsilon_n^{-1}\E \bigl[\bigl\llvert \cA
^{n,m}_{s-}\bigr\rrvert ^2\llvert
W_{s-}-W_{\iota_{n}(s-)}\rrvert ^2 \bigr] \,\dd s
\\
&&\qquad \leq\kappa_1 \int_0^T\E \bigl[
\bigl\llvert \cA^{n,m}_{s-}\bigr\rrvert ^2 \bigr]
\,\dd s\leq\kappa_1 \int_0^T\E
\bigl[\bigl\llvert A^{n}_{s-}\bigr\rrvert ^2
\bigr] \,\dd s
\end{eqnarray*}
for all $n,m\in\N$. The latter integral is uniformly bounded by
Lemma~\ref{boundapp} and the Lipschitz continuity of $a$.

We proceed with the analysis of~(\ref{ggl05}). The expectation in~(\ref{ggl05}) is bounded by twice the sum of
\begin{eqnarray*}
\Sigma_{n,m}^{(1)}&:=&\varepsilon_n^{-1}
\E \biggl[\sup_{t\in
[0,T]}\biggl\llvert \int_0^tD^n_{s-}A^n_{s-}(Y_{s-}-Y_{\iota_{n}(t-)})
\,\dd Y_s
\\
&&\hspace*{57pt} {} -\sigma\int_0^tD^n_{s-}A^{n}_{s-}(W_{s-}-W_{\iota_{n}(s-)})
\,\dd \cY^m_s\biggr\rrvert ^2 \biggr]
\end{eqnarray*}
and
\begin{eqnarray*}
\Sigma_{n,m}^{(2)} &:=& \varepsilon_n^{-1}
\E \biggl[\sup_{t\in
[0,T]}\biggl\llvert \int_0^tD^n_{s-}A^n_{s-}(W_{s-}-W_{\iota_{n}(t-)})
\,\dd \cY^m_s
\\
&&\hspace*{57pt}{}-\int_0^tD^n_{s-}
\cA^{n,m}_{s-}(W_{s-}-W_{\iota_{n}(t-)}) \,\dd
\cY^m_s\biggr\rrvert ^2 \biggr].
\end{eqnarray*}
The term $\Sigma_{n,m}^{(1)}$ is the same as the one appearing
in~(\ref{goodeq2}) when replacing $Y^\varepsilon$ by $\cY^m$. One
can literally translate the proof of~(\ref{goodeq2}) to obtain that
\[
\lim_{m\to\infty}\limsup_{n\to\infty}
\Sigma_{n,m}^{(1)}=0. %
\]
%
By uniform boundedness of $D^n$ and Lemma~\ref{estimateintegral},
there exists a constant $\kappa_2$ not depending on $n,m\in\N$ with
\begin{eqnarray*}
\Sigma_{n,m}^{(2)}&\leq&\kappa_2
\varepsilon_n^{-1}\int_0^T
\E \bigl[\bigl(A^n_{s-}-\cA^{n,m}_{s-}
\bigr)^2(W_{s-}-W_{\iota_{n}(s-)})^2\bigr] \,\dd
s
\\
&\leq&\kappa_2 \int_0^T\E \bigl[
\bigl(A^{n}_{s-}-\cA^{n,m}_{s-}
\bigr)^2 \bigr] \,\dd s
\\
&\leq&2\kappa_2 \int_0^T\E \bigl[
\bigl(A^{n}_{s-}-a(X_{s-})\bigr)^2
\bigr] \,\dd s+2\kappa_2\int_0^T \E
\bigl[\bigl(a(X_{s-})- \cA^{n,m}_{s-}
\bigr)^2 \bigr] \,\dd s,
\end{eqnarray*}
where we have used again that given $\iota_n$ the random variables
$A^n_{s-}-\cA^{n,m}_{s-}$ and $W_{s-}-W_{\iota_{n}(s-)}$ are
independent. The first integral in the previous line tends to zero by
Lipschitz continuity of $a$ and $L^2$-convergence of $\sup_{t\in
[0,T]} \llvert  X^n_t- X_t\rrvert  \to0$ (see Proposition~4.1 of \cite{DerLi14a}).
Further, the second integral satisfies
\[
\limsup_{n\to\infty}\int_0^T \E
\bigl[\bigl(a(X_{s-})- \cA ^{n,m}_{s-}
\bigr)^2 \bigr] \,\dd s \leq\int_0^T
\E\bigl[\1_{[m,\infty
)}\bigl(\llvert X_{s-}\rrvert \bigr)
a(X_{s-})^2 \bigr] \,\dd s %
\]
which tends to zero as $m\to\infty$ since $\sup_{t\in[0,T]} \llvert  X_t\rrvert  $
is square integrable.
\end{pf}

\begin{lemma}
For every $m\in\N$ and $p\geq2$, one has
\[
\sup_{n\in\N}\E \Bigl[\sup_{t\in[0,T]}\bigl\llvert
\cU^{n,n+1,m}_t\bigr\rrvert ^p \Bigr]<\infty.
\]
\end{lemma}

\begin{pf}
Since $\cY^h$ has bounded jumps, it has finite $p$th moment.
$D^{n,n+1}$ is uniformly bounded and by part one of Lemma~\ref
{estimationsde} it suffices to prove that
\[
\E \biggl[ \sup_{t\in[0,T]} \biggl\llvert \varepsilon_n^{-1/2}
\int_0^t D^{n}_{s-}
\cA^{n,m}_{s-}(W_{s-}-W_{\iota_{n}(s-)}) \,\dd
\cY^m_s\biggr\rrvert ^p \biggr]
\]
is uniformly bounded over all $n\in\N$ for fixed $m\in\N$.
Using Lemma~\ref{estimateintegral} and the uniform boundedness of
$D^n$ and $\cA^{n,m}$ over all $n\in\N$, we conclude existence of a
constant $\kappa_3$ such that for every $n\in\N$
\begin{eqnarray*}
&& \E \biggl[ \sup_{t\in[0,T]} \biggl\llvert \int_0^t
D^{n}_{s-}\cA ^{n,m}_{s-}(W_{s-}-W_{\iota_{n}(s-)})
\,\dd\cY^m_s\biggr\rrvert ^p \biggr]
\\
&&\qquad \leq
\kappa_3 \int_0^T\E \bigl[\llvert
W_{s-}-W_{\iota_n(s-)}\rrvert ^p \bigr]\,\dd s\leq
\kappa_3 T \varepsilon_n^{p/2}. %
\end{eqnarray*}\upqed
\end{pf}

\subsection{Proof of the central limit theorems for $X^1,X^2,\ldots$\,}\label{proofCLT1}

In this section we prove Theorems~\ref{theo2402-1} and~\ref{theo2402-2}.
By Proposition~\ref{prop2502-1} and the Lipschitz continuity of $F$
with respect to supremum norm, we conclude that
$(\varepsilon_n^{-1/2} (F(X^{n+1})-F(X^n))\dvtx n\in\N)$ is uniformly
$L^2$-integrable in both settings. In view of the discussion at the
beginning of Section~\ref{proofCLT} it suffices to show that
\[
\lim_{n\to\infty} \Var\bigl(\varepsilon_n^{-1/2}
\bigl(F\bigl(X^{n+1}\bigr)-F\bigl(X^n\bigr)\bigr)\bigr)=\Var\bigl(
\nabla f(AX)\cdot AU\bigr) %
\]
in the first setting and
\[
\lim_{n\to\infty} \Var\bigl(\varepsilon_n^{-1/2}
\bigl(F\bigl(X^{n+1}\bigr)-F\bigl(X^n\bigr)\bigr)\bigr)=\Var
\bigl(f'(X_S)\cdot U_S\bigr) %
\]
in the second setting. By dominated convergence it even suffices to
show weak convergence of the distributions appearing in the variances.
Theorem~\ref{theo2402-1} follows from the following lemma.

\begin{lemma}\label{le2502-6}
Under the assumptions of Theorem~\ref{theo2402-1}, one has
\[
\varepsilon_n^{-1/2} \bigl(F\bigl(X^{n+1}\bigr)-F
\bigl(X^n\bigr)\bigr) \Rightarrow\nabla f(AX)\cdot AU. %
\]
\end{lemma}

\begin{pf}
For $n\in\N$, let $Z_n:=AX^n$ and set $Z=AX$.
Since $Z\in D_f$, almost surely, we conclude that
%
\begin{equation}
\label{eq0602-1} \qquad\lim_{n\to\infty} \varepsilon_n^{-1/2}
\bigl(f(Z_{n})-f(Z)-\nabla f(Z) (Z_{n}-Z)\bigr)=0\qquad
\mbox{in probability}.
\end{equation}
%
Indeed, one has $f(Z_{n})-f(Z)-\nabla f(Z) (Z_{n}-Z)=R_n (Z_n-Z)$ for
appropriate random variables $R_n$ that converge in probability to zero
since $Z_n-Z\to0$, in probability, and $f$ is differentiable in $Z$.
Further, for fixed $\varepsilon>0$ we choose $\delta>0$ large and estimate
\[
\P\bigl(\bigl\llvert \varepsilon_n^{-1/2} R_n
(Z_n-Z)\bigr\rrvert >\varepsilon\bigr) \leq\P \bigl(\llvert
R_n\rrvert >\varepsilon/\delta\bigr)+\P\bigl(\bigl\llvert
\varepsilon_n^{-1/2}(Z_n-Z)\bigr\rrvert >\delta
\bigr). %
\]
The first summand converges to zero as $n\to\infty$ and the second
term can be made uniformly arbitrarily small over $n$ by choosing
$\delta$ sufficiently large due to tightness of the sequence
$(\varepsilon_n^{-1/2} (Z_n-Z))_{n\in\N}$.
Equation~(\ref{eq0602-1}) remains true when replacing $Z_n$ by
$Z_{n+1}$ and we conclude that
\[
\lim_{n\to\infty} \varepsilon_n^{-1/2}
\bigl(f(Z_{n+1})- f(Z_n) -\nabla f(Z) (Z_{n+1}-Z_n)\bigr)=0
\qquad\mbox{in probability}. %
\]
By Theorem~\ref{twosuccessive} and the fact that $A$ is continuous in
$\P_U$-almost every point, we conclude that
\[
\bigl(Y,A\varepsilon_n^{-1/2} \bigl(X^{n+1}-X^{n}
\bigr)\bigr)\Rightarrow(Y,AU)
\]
and, hence,
\[
\varepsilon_n^{-1/2}
(Z_{n+1}-Z_n)\stackrel {\operatorname{stably}} {
\Longrightarrow}AU, %
\]
by Lemma~\ref{stablethm1}.
Consequently, since $\nabla f(Z)$ is $\sigma(Y)$-measurable we get
\[
\bigl(\nabla f(Z), \varepsilon_n^{-1/2}
(Z_{n+1}-Z_n)\bigr)\Rightarrow\bigl(\nabla f(Z),AU\bigr)
\]
and the proof is completed by noticing that the scalar product is continuous.
\end{pf}

Analogously, Theorem~\ref{theo2402-2} is a consequence of the
following lemma.

\begin{lemma}Under the assumptions of Theorem~\ref{theo2402-2}, one has
\[
\varepsilon_n^{-1/2} \Bigl(f \Bigl(\sup
_{s\in[0,T]} X^{n+1}_s \Bigr)-f \Bigl(\sup
_{s\in[0,T]} X_s^n \Bigr) \Bigr)
\Rightarrow f'(X_S)\cdot U_S,
\]
where $S$ denotes the time where $X$ attains its maximum.
\end{lemma}

\begin{pf}
By Lemma~\ref{le2502-3}, there exists s unique time $S$ at which $X$
attains its maximum and by Lemma~\ref{le2502-4} one has
\[
\varepsilon_n^{-1/2} \Bigl(\sup_{s\in[0,T]}
X_s^{n+1}-\sup_{s\in
[0,T]}
X_s^{n} \Bigr) - U^{n,n+1}_S \to0
\qquad\mbox{in probability}. %
\]
By Theorem~\ref{twosuccessive} and Lemma~\ref{stablethm1}, one has
\[
\bigl(Y, S,U^{n,n+1}\bigr) \Rightarrow(Y, S,U) %
\]
and the function $[0,T]\times\ID(\R)\to\R, (s,u)\mapsto u_s$ is
continuous in $\P_{S,U}$-almost all $(s,u)$ since $U$ is almost surely
continuous in $S$ by Lemma~\ref{le2502-3}. Consequently,
\[
\bigl(Y, U^{n,n+1}_S\bigr)\Rightarrow(Y, U_S)
\]
and, hence,
\[
\varepsilon_n^{-1/2} \Bigl(\sup
_{s\in[0,T]} X_s^{n+1}-\sup
_{s\in
[0,T]} X_s^{n} \Bigr)\stackrel{
\operatorname{stably}} {\Longrightarrow}U_S. %
\]
The rest follows as in the proof of Lemma~\ref{le2502-6}.
\end{pf}


\subsection{Proof of Theorem~\texorpdfstring{\protect\ref{theo0303-1}}{1.11}}\label{sec43}

1\textit{st step}. Denote by $\cE=(\cE_t)_{t\in[0,T]}$ the stochastic
exponential of $(\int_0^t a'(X_{s-}) \,\dd Y_s)_{t\in[0,T]}$. In
particular, $\cE$ does not hit zero with probability one; see, for
instance, \cite{jacodshibook}, Theorem~1.4.61.
In the first step, we show that $\E[U_sU_t\mid Y]=\Upsilon^2  \phi
_{s,t}(Y)$, where
%
\begin{eqnarray}\label{eq2802-4}
\phi_{s,t}(Y) &=& \sigma^4 \cE_s
\cE_t \int_0^s\frac
{(aa')(X_{u-})^2}{\cE_{u-}^2}
\,\dd u
\nonumber\\[-8pt]\\[-8pt]\nonumber
&&{} + \sigma^2 \lim_{\delta
\downarrow0} \cE_s
\cE_t\mathop{\sum_{u\in(0,s]:}}_{\llvert  \Delta
Y_u\rrvert  \geq\delta}
\frac{(aa')(X_{u-})^2 \Delta Y_u^2}{(1+a'(X_{u-})
\Delta Y_u)^2\cE_{u-}^2}
\end{eqnarray}
and the limit is taken in ucp.

We define $\overline L=(\overline L_t)_{t\in[0,T]}$ by
\[
\overline L_t=\sig^2 \Upsilon B_t+ \lim
_{\delta\downarrow0} \mathop {\sum_{s\in(0,t]:}}_{\llvert  \Delta Y_s\rrvert  \geq\delta}
\frac{1}{1+a'(X_{s-})
\Delta Y_s} \Delta L_s %
\]
and note that the process is well-defined since the denominator does
not attain the value zero by assumption.
Using the product rule and independence of $W$ and $B$, it is straight
forward to verify that
\[
\biggl(\cE_t \int_0^t
\frac{(aa')(X_{s-})}{\cE_{s-}} \,\dd\overline L_s \biggr)_{t\in[0,T]} %
\]
solves the stochastic integral equation~(\ref{limitproce}) and by
strong uniqueness of the solution equals~$U$, almost surely.
We write
\[
U_t= \sigma^2 \Upsilon\underbrace{\cE_t \int
_0^t \frac
{(aa')(X_{s-})}{\cE_{s-}} \,\dd
B_s}_{=:Z_t}+ \lim_{\delta\downarrow
0} \underbrace{
\cE_t\mathop{\sum_{s\in(0,t]:}}_{\llvert  \Delta
Y_s\rrvert  \geq\delta}
\frac{(aa')(X_{s-})}{(1+a'(X_{s-}) \Delta Y_s)\cE
_{s-}} \Delta L_s}_{Z^{(\delta)}_t}
\]
and note that given $Y$ the processes $Z$ and $Z^{(\delta)}$ are
independent and have expectation zero.
Further, for $0\leq s\leq t\leq T$ one has
\[
\E[Z_sZ_t\mid Y]= \cE_s \cE_t
\int_0^s\frac{(aa')(X_{u-})^2}{\cE
_{u-}^2} \,\dd u %
\]
and
\[
\E\bigl[Z^{(\delta)}_sZ^{(\delta)}_t\mid Y
\bigr]=\cE_s\cE_t\mathop{\sum
_{u\in
(0,s]:}}_{\llvert  \Delta Y_u\rrvert  \geq\delta} \frac{(aa')(X_{u-})^2 \Delta
Y_u^2}{(1+a'(X_{u-}) \Delta Y_u)^2\cE_{u-}^2} \E\bigl[
\sigma^2_u\bigr]. %
\]
One easily computes that $\E[\sigma^2_u]=\sigma^2 \Upsilon^2$.
Altogether, it follows the wanted statement.

2\textit{nd step}. Let $A=(A_1,\ldots,A_d)\dvtx \ID(\R)\to\R^d$ be a
linear map of integral type meaning that there are finite signed
measures $\mu_1,\ldots,\mu_d$ on $[0,T]$ with
\[
A_jx=\int_0^Tx_s
\,\dd\mu_j(s). %
\]
Then by conditional Fubini and step one,
\begin{eqnarray*}
&& \Var\bigl[\nabla f(AX) \cdot AU\bigr]
\\
&&\qquad = \sum_{i,j=1}^d
\E\bigl[\partial_i f(AX) A_iU \,\partial_j
f(AX) A_jU\bigr]
\\
&&\qquad =\sum_{i,j=1}^d \E \biggl[
\partial_i f(AX) \,\partial_j f(AX) \E \biggl[ \int
_{[0,T]^2} U_u U_v \,\dd
\mu_i\otimes\mu_j(u,v)\Big| Y \biggr] \biggr]
\\
&&\qquad =\Upsilon^2 \sum_{i,j=1}^d \E
\biggl[\partial_i f(AX) \,\partial_j f(AX) \int
_{[0,T]^2} \phi_{u,v}(Y) \,\dd\mu_i\otimes
\mu _j(u,v) \biggr].
\end{eqnarray*}

3\textit{rd step}. The supremum dependent case follows by noticing that
step one remains valid when choosing $s=t=S$ since $S$ is $\sigma
(Y)$-measurable.



\setcounter{theorem}{0}
\setcounter{equation}{0}
\begin{appendix}\label{append}
\section*{Appendix}

\subsection{Stable and weak convergence}\label{sec1803-01}

We briefly introduce the concept of stable convergence first appearing
in R\'enyi \cite{ren}. 

\begin{definition}Let $\cF^0$ denote a sub-$\sigma$-field of $\cF$.
A sequence $(Z_n)_{n\in\N}$ of $\cF^0$-measurable random variables
taking values in a Polish space $E$ converges stably with respect to
$\cF^0$ to an $E$-valued $\cF$-measurable random variable $Z$,
if for every $A\in\cF^0$ and continuous and bounded function $f\dvtx E\to
\R$
\[
\lim_{n\to\infty} \E\bigl[\1_A f(Z_n)
\bigr] =\E\bigl[\1_A f(Z)\bigr]. %
\]
We briefly write $Z^n\stackrel{\operatorname{stably}}{\Longrightarrow}Z$.
\end{definition}

Stable convergence admits various equivalent definitions.

\begin{theorem}\label{stablethm1}
Let $(Z_n)$ and $Z$ be $\cF^0$-measurable, respectively, $\cF
$-measu\-rable, random variables taking values in a Polish space $E$.
The following statements are equivalent:
\begin{longlist}[2.]
\item[1.] $Z_n\stackrel{\operatorname{stably}}{\Longrightarrow}Z$ with
respect to $\cF^0$,\vspace*{1pt}

\item[2.] for all bounded $\cF^0$-measurable random variables $U$ and all
bounded and continuous functions $f\dvtx E\to\R$ one has
%
\begin{equation}
\label{eq1001-1} \lim_{n\to\infty} \E\bigl[U f(Z_n)\bigr]
=\E\bigl[U f(Z)\bigr].
\end{equation}
\end{longlist}
If $\cF^0=\sig(Y)$ for a random variable $Y$ taking values in a
Polish space $E'$, then stable convergence is equivalent to weak convergence
%
\begin{equation}
\label{eq2001-2} (Y,Z_n)\Rightarrow(Y,Z)\qquad\mbox{in } E\times
E'.
\end{equation}
%
\end{theorem}

\begin{pf} The first equivalence is an immediate consequence of the
fact that the set of $\cF^0$-measurable random variables $U$ for which
\[
\lim_{n\to\infty} \E\bigl[U f(Z_n)\bigr] =\E\bigl[U
f(Z)\bigr] %
\]
is true is linear and closed with respect to $L^1$-norm.
Further,~(\ref{eq2001-2}) implies $Z_n\stackrel{\operatorname
{stably}}{\Longrightarrow}Z$ since the $L^1$-closure of random
variables $g(Y)$ with $g\dvtx E'\to\R$ bounded and continuous contains all
indicators $\1_A$ with $A\in\cF^0$. Conversely, assuming
$Z_n\stackrel{\operatorname{stably}}{\Longrightarrow}Z$, the
sequence of random variables $((Y,Z_n)\dvtx n\in\N)$ is tight in the
product topology and for any $g\dvtx E'\to\R$ bounded and continuous one has
$\E[g(Y)f(Z_n) ]\to\E[g(Y)f(Z)]$ which implies that
$(Y,Z_n)\Rightarrow(Y,Z)$. The last statement is proved in complete
analogy with the proof of the corresponding statement for weak convergence.
\end{pf}

As the latter theorem shows, stable and weak convergence are intimately
connected and we will make use of results of Jacod and Protter~\cite
{jacodprotter} on weak convergence for stochastic differential
equations. For the statement, we need the concept of uniform tightness.

\begin{definition} Let $(\cF_t)_{t\in[0,T]}$ be a filtration and
$(Z^n\dvtx n\in\N)$ be a sequence of c\`adl\`ag $(\cF_t)$-semimartingales.
For $\delta>0$ we represent each semi-martingale uniquely in the form
\[
Z^{n}_t=Z^{n}_0+A^{n,\delta}_t+M^{n,\delta}_t+
\sum_{s\leq t}\Delta Z^{n}_s
\1_{\{\llvert  \Delta Z^{n}_s\rrvert  >\delta\}}\qquad\mbox{for } t\in[0,T], %
\]
where $A^{n,\delta}=(A^{n,\delta}_t)_{t\in[0,T]}$ is a c\`adl\`ag
predictable process of finite variation and $M=(M^{n,\delta}_t)_{t\in
[0,T]}$ is a c\`adl\`ag local martingale, both processes starting in
zero. We say that $(Z^{n}\dvtx n\in\N)$ is \emph{uniformly tight}, if the
sequence,
\[
\bigl\langle M^{n,\delta},M^{n,\delta}\bigr\rangle_T+\int
_0^T\bigl\llvert \,\dd A^{n,\delta
}_s
\bigr\rrvert +\sum_{0\leq s\leq T}\bigl\llvert \Delta
Z^{n,i}\bigr\rrvert \1_{\{\llvert  \Delta
Z^{n,i}_s\rrvert  >\delta\}} %
\]
is tight. The definition does not depend on the particular choice of
$\delta$. Multivariate processes are called uniformly tight if each
component is uniformly tight.
\end{definition}

We cite \cite{jacodprotter}, Theorem 2.3, which is a consequence
of~\cite{kurzprotter}.

\begin{theorem}\label{weakproperty1} Let $Z,Z^1,Z^2,\ldots$ be
c\`adl\`ag one-dimensional semimartingales and $H$ be a c\`adl\`ag
one-dimensional adapted process. If:
\begin{longlist}[(ii)]
\item[(i)] $(Z^{n}\dvtx n\in\N)$ is uniformly tight and
\item[(ii)] $((H,Z^n)\dvtx n\in\N)\Rightarrow(H,Z)$ in $D(\R^{2})$,
\end{longlist}
then
\[
\biggl(H,Z^n,\int_0^{\cdot}
H_{s-}\,\dd Z^n_s\dvtx n\in\N \biggr)\Rightarrow
\biggl(H,Z,\int_0^{\cdot} H_{s-}\,\dd
Z_s \biggr)\qquad \mbox{in } D\bigl(\R^{3}\bigr).
\]
\end{theorem}

We state a consequence of \cite{kurzprotter}, Theorem 8.2.

\begin{theorem}\label{weakproperty2}
Let $H, Z, Z^1,Z^2, \ldots$ be as in the previous theorem. Further, let
$Y$ be an adapted c\`adl\`ag semimartingale. We define
$U^n:=(U^n_t)_{t\in[0,T]}$ and $U:=(U_t)_{t\in[0,T]}$ by
\[
U^n_t=Z^n_t+\int
_0^tU^n_{s-}H_{s-}
\,\dd Y_s,\qquad U_t=Z_t+\int
_0^tU_{s-}H_{s-}\,\dd
Y_s\qquad\mbox{for }t\in[0,T]. %
\]
If
\[
\biggl(Z^n,\int_0^{\cdot}H_{s-}
\,\dd Y_s \biggr)\Rightarrow \biggl(Z,\int_0^{\cdot}H_{s-}
\,\dd Y_s \biggr)\qquad\mbox{in }D\bigl(\R^{2}\bigr),
\]
then
\[
\biggl(Z^n,\int_0^{\cdot}H_{s-}
\,\dd Y_s, U^n \biggr)\Rightarrow \biggl(Z,\int
_0^{\cdot}H_{s-}\,\dd Y_s,U
\biggr)\qquad\mbox{in }D\bigl(\R^{3}\bigr). %
\]
\end{theorem}

The definition of uniform tightness and the two theorems above have
natural extension to the multivariate setting and we refer the reader
to \cite{kurzprotter} for more details. Further results about stable
convergence of stochastic process can be found in \cite{jacod2} and
\cite{jacodshibook}.

A helpful lemma in the treatment of weak convergence is the following.

\begin{lemma}\label{le2901-1} Let $A,A^1,A^2,\ldots$ be processes with
trajectories in $\ID(\R^d)$.
\begin{enumerate}
\item Suppose that for every $m\in\N$, $A^m,A^{1,m},A^{2,m},\ldots$
are processes with trajectories in $\ID(\R^d)$ such that:
\begin{longlist}[(a)]
\item[(a)]$\forall\delta>0$:  $\lim_{m\to\infty}\limsup_{n\to\infty}
\P(\sup_{t\in[0,T]}\llvert  A_t^{n,m}-A_t^n\rrvert  >\delta)=0$,\vspace*{1pt}
\item[(b)]$\lim_{m\to\infty} \P(\sup_{t\in[0,T]}\llvert  A_t^{m}-A_t\rrvert  >\delta)=0$.
\end{longlist}
Provided that one has convergence $A^{n,m}\Rightarrow A^m$ for every
$m\in\N$, it is also true that
\[
A^n\Rightarrow A. %
\]
\item Suppose that $B^1,B^2,\ldots$ are processes with trajectories in
$\ID(\R^d)$ such that for all $\delta>0$
\[
\lim_{n\to\infty} \P\Bigl(\sup_{t\in[0,T]}\bigl\llvert
B^n_t-A^n_t\bigr\rrvert >\delta
\Bigr)=0. %
\]
Then one has weak convergence $A^n\Rightarrow A$ if and only if
$B^n\Rightarrow A$.
\end{enumerate}
\end{lemma}

\begin{pf}
To prove weak convergence on $\ID(\R^d)$ it suffices to consider
bounded and continuous test functions $f\dvtx \ID(\R^d)\to\R$ that are
additionally Lipschitz continuous with respect to supremum norm. Using
this, it is elementary to verify the first statement. Further, the
second statement is an immediate consequence of the first one.
\end{pf}

\begin{remark}\label{re2901-1}In general, we call approximations
$A^m,A^{1,m},A^{2,m},\ldots$ with properties (a) and (b) of part one of
the lemma \emph{good approximations} for $A,A^1,A^2,\ldots.$ Further,
approximations $B^1,B^2,\ldots$ as in part two will be called \emph
{asymptotically equivalent in ucp} to $A^1,A^2,\ldots.$
\end{remark}

\subsection{Auxiliary estimates}

We will make use of the following analogue of Lemma~\ref{le2901-1} for
tightness.

\begin{lemma}\label{lemmauniforminteg}
Let $(A_n)_{n\in\N}$ and, for every $m\in\N$, $(A^{(m)}_n)_{n\in\N
}$ be sequences of $L^2$-integrable random variables.
If
\[
\lim_{m\to\infty}\limsup_{n\to\infty}\E\bigl[\bigl
\llvert A_n-A^{(m)}_n\bigr\rrvert
^2\bigr]=0 %
\]
and, for every $m\in\N$, the sequence $(A^{(m)}_n)_{n\in\N}$ is
uniformly $\L^2$-integrable, then also the sequence $(A_n)_{n\in\N}$
is uniformly $\L^2$-integrable. In particular, if there is a sequence
$(B_n)_{n\in\N}$ of uniformly $\L^2$-integrable random variables with
\[
\lim_{n\to\infty}\E\bigl[\llvert B_n-A_n
\rrvert ^2\bigr]=0, %
\]
then $(A_n)_{n\in\N}$ is uniformly $\L^2$-integrable.
\end{lemma}

\begin{pf}
For $\eta>0$ and $n,m\in\N$, one has
\begin{eqnarray*}
\E \bigl[\llvert A_n\rrvert ^2\1_{\{\llvert  A_n\rrvert  \geq\eta\}}
\bigr]&\leq&2 \E \bigl[\bigl\llvert A_n-A_n^{(m)}
\bigr\rrvert ^2 \bigr] + 2 \E \bigl[\bigl\llvert A_n^{(m)}
\bigr\rrvert ^2 \1_{\{\llvert  A_n\rrvert  \geq
\eta\}} \bigr]
\\
&\leq&2 \E \bigl[\bigl\llvert A_n-A_n^{(m)}
\bigr\rrvert ^2 \bigr] + 2 \E \bigl[\bigl\llvert A_n^{(m)}
\bigr\rrvert ^2 \1_{\{\llvert  A^{(m)}_n\rrvert  \geq\eta/2\}} \bigr]
\\
&&{} + 2 \E \bigl[\bigl\llvert A_n^{(m)}\bigr\rrvert
^2 \1_{\{ \llvert  A_n^{(m)}\rrvert  <\eta/2,
\llvert  A_n-A_n^{(m)}\rrvert  \geq\eta/2\}} \bigr]
\\
&\leq&2 \E \bigl[\bigl\llvert A_n-A_n^{(m)}
\bigr\rrvert ^2 \bigr] + 2 \E \bigl[\bigl\llvert A_n^{(m)}
\bigr\rrvert ^2 \1_{\{\llvert  A^{(m)}_n\rrvert  \geq\eta/2\}} \bigr]
\\
&&{}+\frac{\eta
^2}2 \P\bigl(
\bigl\llvert A_n-A_n^{(m)}\bigr\rrvert \geq
\eta/2\bigr)
\\
&\leq&4 \E \bigl[\bigl\llvert A_n-A_n^{(m)}
\bigr\rrvert ^2 \bigr] + 2 \E \bigl[\bigl\llvert A_n^{(m)}
\bigr\rrvert ^2 \1_{\{\llvert  A^{(m)}_n\rrvert  \geq\eta/2\}} \bigr],
\end{eqnarray*}
where we used Chebychew's inequality in the last step.
Let now $\varepsilon>0$. By assumption, we can choose $m$ sufficiently
large such that for all large $n$, say for $n\geq n_0$, $4 \E
[\llvert  A_n-A_n^{(m)}\rrvert  ^2 ]\leq\varepsilon/2$. Further,\vspace*{1pt} by the uniform
$L^2$-integrability of $(A_n^{(m)})_{n\in\N}$ we can choose $\eta$
large to ensure that for all $n\in\N$, $2  \E [\llvert  A_n^{(m)}\rrvert  ^2 \1
_{\{\llvert  A^{(m)}_n\rrvert  \geq\eta/2\}} ]\leq\varepsilon/2$ so that $\E
 [\llvert  A_n\rrvert  ^2\1_{\{\llvert  A_n\rrvert  \geq\eta\}} ]\leq\varepsilon$ for
$n\geq n_0$. For $n=1,\ldots,n_0-1$ this estimate remains true for a
sufficiently enlarged $\eta$, since\vspace*{1pt} finitely many $L^2$-integrable
random variables are always uniformly $L^2$-integrable.
\end{pf}

\begin{lemma}\label{le1702-4}Let $A_1,A_2,\ldots$ be real random
variables and let $(\varepsilon_k)_{k\in\N}$ satisfy \textup{(ML1)} (see
Section~\ref{secMLintro}), and $L(\delta)$ and $n_k(\delta)$ be as
in~(\ref{eq1103-2}). Suppose that:
\begin{longlist}[2.]
\item[1.]$\Var(\varepsilon_{k-1}^{-1/2} A_k)\to\zeta$ and
\item[2.]$(\varepsilon_{k-1}^{-1/2} A_k\dvtx k\in\N)$ is $L^2$-uniformly integrable.
\end{longlist}
Denote by $(A_{k,j}\dvtx k,j\in\N)$ independent random variables with $\cL
(A_{k,j})=\cL(A_k)$. The random variables $(\widehat S_\delta\dvtx \delta
\in(0,1))$ given by
\[
\widehat S_\delta:= \sum_{k=1}^{L(\delta)}
\frac{1}{n_k(\delta)} \sum_{j=1}^{n_k(\delta)}
A_{k,j} %
\]
satisfy
\[
\delta^{-1}\bigl(\widehat S_\delta-\E[\widehat S_\delta]\bigr) \Rightarrow \cN(0,\zeta). %
\]
\end{lemma}

\begin{pf} Without loss of generality, we can and will assume that the
random variables $A_1,A_2,\ldots$ have zero mean.

\begin{longlist}[1\textit{st step}.]
\item[1\textit{st step}.] We first show that the variance of $\widehat S_\delta
$ converges. One has
\[
\Var(\widehat S_\delta)=\sum_{k=1}^{L(\delta)}
\frac{1}{n_k(\delta
)} \Var(A_k)=\sum_{k=1}^{L(\delta)}
\underbrace{ \biggl\lfloor\frac
{\delta^2}{L(\delta) \varepsilon_{k-1} } \biggr\rfloor\varepsilon
_{k-1}}_{=:a_{k,\delta}}\frac{\Var(A_k)}{\varepsilon_{k-1}}.
\]
It is elementary to verify that $ \sum_{k=1}^{L(\delta)}( a_{k,\delta
} \delta^{-2}-L(\delta)^{-1})\to0$ as $\delta\downarrow0$. By the
boundedness of $({\Var(A_k)}/{\varepsilon_{k-1}})_{k\in\N}$ one has
\[
\Biggl\llvert \delta^{-2} \Var(\widehat S_\delta)-
\frac{1}{L(\delta)}\sum_{k=1}^{L(\delta)}
\frac{\Var(A_k)}{\varepsilon_{k-1}}\Biggr\rrvert \to0 %
\]
and we get that $\lim_{\delta\downarrow0} \Var(\delta^{-1}\widehat S_\delta)=\zeta$ since the C\'esaro mean of a convergent sequence
converges to its limit.

\item[2\textit{nd step}.] In view of the Lindeberg condition (see, e.g., \cite{Kallenberg2nd}, Theorem 5.12), it suffices to verify
that for arbitrarily fixed $\kappa>0$ one has
\[
\Sigma(\delta):=\sum_{k=1}^{L(\delta)} \sum
_{j=1}^{n_k^{(\delta
)}} \E \biggl[ \biggl(
\frac{A_{k,j}}{\delta n_k^{(\delta)}} \biggr)^2 \1_{\{\llvert  A_{k,j}/(\delta n_k^{(\delta)})\rrvert  >\kappa\}} \biggr]\to0 \qquad
\mbox{as } \delta\downarrow0. %
\]
We estimate
\[
\Sigma(\delta) \leq\delta^{-2} \sum_{k=1}^{L(\delta)}
\frac
{\varepsilon_{k-1}}{n_k^{(\delta)} } \E \biggl[ \frac
{A_{k}^2}{\varepsilon_{k-1}} \1_{\{{\llvert  A_{k}\rrvert  }/{\sqrt{\varepsilon
_{k-1}}}>{\kappa\delta n_k^{(\delta)}}/{\sqrt{\varepsilon
_{k-1}}} \}} \biggr]
\]
and note that for $k=1,\ldots,L(\delta)$
\[
\varepsilon_{k-1}\geq\varepsilon_{L(\delta)-1} = T M^{-L(\delta
)+1}
\geq T \delta^2, %
\]
where we used that $\alpha\geq1/2$ in the previous step. Hence, for
these $k$, one has $\delta n_k^{(\delta)}/\sqrt{\varepsilon
_{k-1}}\geq\delta^{-1} \sqrt{\varepsilon_{k-1}} L(\delta) \geq
\sqrt T L(\delta)$. Consequently,
\[
\Sigma(\delta) \leq\delta^{-2} \sum_{k=1}^{L(\delta)}
\frac
{\varepsilon_{k-1}}{n_k^{(\delta)} } \E \biggl[ \frac
{A_{k}^2}{\varepsilon_{k-1}} \1_{\{{\llvert  A_{k}\rrvert  }/{\sqrt{\varepsilon
_{k-1}}}>\kappa\sqrt T L(\delta) \}} \biggr].
\]
By uniform $L^2$-integrability of $(A_k/\sqrt{\varepsilon
_{k-1}})_{k\in\N}$ and the fact that $L(\delta)\to\infty$, we get that
\[
\E \biggl[ \frac{A_{k}^2}{\varepsilon_{k-1}} \1_{\{
{\llvert  A_{k}\rrvert  }/{\sqrt{\varepsilon_{k-1}}}>\kappa\sqrt T L(\delta) \}
} \biggr]\leq a(\delta)\qquad
\mbox{for } k=1,\ldots,L(\delta), %
\]
with $(a_\delta)_{\delta\in(0,1)}$ being positive reals with $\lim_{\delta\downarrow0} a_\delta=0$. Hence, $\Sigma(\delta)\leq\break
a_\delta\delta^{-2} \sum_{k=1}^{L(\delta)} \frac{\varepsilon
_{k-1}}{n_k^{(\delta)} }$ and we remark that the analysis of step one
yields equally well that $\delta^{-2} \sum_{k=1}^{L(\delta)} \frac
{\varepsilon_{k-1}}{n_k^{(\delta)} }$ converges to a finite limit.\quad\qed
\end{longlist}\noqed
\end{pf}

\subsection{Estimates for L\'evy-driven SDEs}

Let $Y=(Y_t)_{t\in[0,T]}$ denote a square integrable L\'evy process
with triplet $(b,\sig^2,\nu)$.

\begin{lemma}\label{driftdisappear}
Let $(\varepsilon_n)$ and $(h_n)$ be positive decreasing sequences
such that
\[
\sup_{n\in\N} \nu\bigl(B(0,h_n)^c\bigr)
\varepsilon_n <\infty. %
\]
One has
%
\begin{equation}
\label{eq1902-1} \varepsilon_n \biggl(\int_{B(0,h_n)^c}x
\nu(\dd x) \biggr)^2\rightarrow0\qquad\mbox{as }n\rightarrow\infty.
\end{equation}
Further, if the limit $\lim_{n\to\infty} \nu(B(0,h_n))^c
\varepsilon_n=:\theta$ exists and is strictly positive, then
$\lim_{n\to\infty} h_n/\sqrt{\varepsilon_n}=0$. If additionally
$\int x^2 \log^2(1+1/x) \nu(\dd x)<\infty$, then
%
\begin{eqnarray}\label{eq1703-1}
\lim_{n\to\infty} \int_{B(0,h_n)}
x^2 \nu(\dd x) \log^2 \biggl(1+\frac{1}{\varepsilon_n}
\biggr)&=&0\quad\mbox{and}
\nonumber\\[-8pt]\\[-8pt]\nonumber
\lim_{n\to\infty} \frac{h_n^2}{\varepsilon_n}
\log^2 \biggl(1+\frac{1}{\varepsilon_n} \biggr) &=&0.
\end{eqnarray}
\end{lemma}

\begin{pf}
One has for fixed $h>0$ for all $n\in\N$ that
\begin{eqnarray*}
&& \varepsilon_n \biggl(\int_{B(0,h_n)^c} x \nu(\dd x)
\biggr)^2
\\
&&\qquad \leq 2 \varepsilon_n \biggl(\int
_{B(0,h)^c} x \nu(\dd x) \biggr)^2+2
\varepsilon_n \biggl(\int_{B(0,h)\setminus B(0,h_n)} x \nu (\dd x)
\biggr)^2. %
\end{eqnarray*}
The first term on the right-hand side tends to zero since $\varepsilon
_n$ tends to zero. Further, the Cauchy--Schwarz inequality yields for
the second term
\[
\varepsilon_n \biggl(\int_{B(0,h)\setminus B(0,h_n)}x \nu(\dd x)
\biggr)^2\leq\varepsilon_n \nu\bigl(B(0,h_n)^c
\bigr) \int_{B(0, h)} x^2 \nu(\dd x). %
\]
By assumption, $(\varepsilon_n \nu(B(0,h_n)^c))$ is uniformly bounded
and by choosing $h$ arbitrarily small we can make the integral as small
as we wish. This proves~(\ref{eq1902-1}).

We assume that $\lim_{n\to\infty} \nu(B(0,h_n))^c \varepsilon
_n=:\theta>0$. The second statement follows by noting that
\[
\theta^2 \frac{h_n^2}{\varepsilon_n} \sim\varepsilon_n
h_n^2 \nu\bigl(B(0,h_n)^c
\bigr)^2 \leq\varepsilon_n \biggl(\int
_{B(0,h_n)^c}x \nu (\dd x) \biggr)^2\to0. %
\]
The first estimate in~(\ref{eq1703-1}) follows from
\[
\int_{B(0,h_n)}x^2\nu(\dd x) \leq\underbrace{\int
_{B(0,h_n)}x^2\log^2 \biggl(1+
\frac{1}x \biggr) \nu(\dd x)}_{\to0} \bigl({\log(1+1/h_n)}
\bigr)^{-2} %
\]
and recalling that $h_n/\sqrt{\varepsilon_n}\to0$. The second
estimate in~(\ref{eq1703-1}) follows in complete analogy to the proof
of~(\ref{eq1902-1}).
\end{pf}

\begin{lemma}\label{estimateintegral}
Let $p\geq2$ and suppose that $\E[\llvert  Y_T\rrvert  ^{p}]<\infty$. Then there
exists a finite constant $\kappa$ such that for every predictable
process $H$ one has
\[
\E \biggl[\sup_{t\in[0,T]}\biggl\llvert \int_0^t
H_s \,\dd Y_s\biggr\rrvert ^{p} \biggr]\leq
\kappa\int_0^T\E\bigl[\llvert H_s
\rrvert ^{p}\bigr] \,\dd s.
\]
If $p=2$, one can choose $\kappa=2b^2T+8(\sigma^2+\int x^2 \nu(\dd x))$.
\end{lemma}

\begin{pf}
The proof is standard; see, for instance, \cite{protterphilip}, Theorem~V.66. The explicit constant in the $p=2$ case can be
deduced with Doob's $L^2$-inequality and the Cauchy--Schwarz inequality.
\end{pf}

\begin{lemma}\label{boundapp} Irrespective of the choice of the
parameters $(\varepsilon_n)$ and $(h_n)$, one has
\[
\sup_{n\in\N} \E \Bigl[\sup_{t\in[0,T]}\bigl\llvert
X^n_t\bigr\rrvert ^{2} \Bigr]<\infty.
\]
\end{lemma}

The proof of the lemma is standard and can be found, for instance,
in~\cite{kohatsutankov}, Lemma 8. 
%
%

\subsection{Perturbation estimates for SDEs}\label{secperturb}

In this section, we collect perturbation estimates for solutions of
stochastic differential equations.
For $n,m\in\N$, we denote by $\cZ^n$, $\overline\cZ^{n}$, $\cZ^{n,m}$
and $\overline\cZ^{n,m}$ c\`adl\`ag semimartingales and by~$Y$ a square
integrable L\'evy process all with respect to the same filtration.
Further, let $H^{n}$, $H^{n,m}$ and $H$ be c\`agl\`ad adapted processes.
We represent $Y$ as in~(\ref{firstdecom}) and consider as
approximations the processes $Y^m=(Y^m_t)_{t\in[0,T]}$ given by
\[
Y^m_t=bt+\sigma W_t+ \lim
_{\delta\downarrow0} \int_{(0,t]\times
(V_m\setminus B(0,\delta))} x \,\dd\overline\Pi(s,x),
\]
where $V_1,V_2,\ldots$ denote an increasing sequence of Borel sets with
$\bigcup_{m\in\N} V_m=\R\setminus\{0\}$.

In the first part of the subsection, we derive perturbation estimates
for the processes
$\cU^{n,m}=(\cU^{n,m}_t)_{t\in[0,T]}$ and $\overline\cU^{n,m}=(\overline\cU
^{n,m}_t)_{t\in[0,T]}$ given as solutions to
\[
\cU^{n,m}_t=\int_0^t
\cU^{n,m}_{s-}H^{n,m}_{s} \,\dd
Y^m_s+\cZ^{n,m}_t %
\]
and
\[
\overline\cU^{n,m}_t=\int_0^t
\overline\cU^{n,m}_{s-}H^{n,m}_{s} \,\dd
Y_s+\overline\cZ^{n,m}_t. %
\]

\begin{lemma}\label{estimatesde2}
Suppose that
%
\begin{equation}
\label{eq1802-1}\sup_{t\in[0,T]}\bigl\llvert H^{n,m}_t
\bigr\rrvert \quad \mbox{and}\quad\E \Bigl[\sup_{t\in[0,T]}\bigl
\llvert \cZ^{n,m}_t\bigr\rrvert ^2 \Bigr]
\end{equation}
are uniformly bounded over all $n,m\in\N$. Then
\[
\sup_{n,m\in\N} \E \Bigl[\sup_{t\in[0,T]}\bigl\llvert
\cU^{n,m}_t\bigr\rrvert ^2 \Bigr]<\infty.
\]
\end{lemma}

\begin{pf}Suppose that the expressions in~(\ref{eq1802-1}) are bounded
by~$\kappa_1$, denote by~$\cT$ a stopping time\vspace*{1pt} and define $z_\cT
(t)=\E[\sup_{s\in[0,t\wedge\cT]}\llvert  \cU^{n,m}_s\rrvert  ^2]$ for $t\in
[0,T]$. By Lemma~\ref{estimateintegral}, there exists a finite
constant~$\kappa_2$ such that
\begin{eqnarray*}
z_\cT(t)&\leq&2\kappa_2 \int_0^{t}
\E \bigl[\1_{\{s\leq\cT\}} \bigl\llvert \overline\cU^{n,m}_{s-}
\bigr\rrvert ^2\bigl\llvert H^{n,m}_s\bigr
\rrvert ^2 \bigr] \,\dd s+2 \E \Bigl[\sup_{s\in[0,t]}\bigl
\llvert \cZ^{n,m}_s\bigr\rrvert ^2 \Bigr]
\\
&\leq&2\kappa_2\kappa^2_1 \int
_0^t z_\cT(s) \,\dd s+2
\kappa_1.
\end{eqnarray*}
We replace $\cT$ by a localising sequence $(\cT_k)_{k\in\N}$ of
stopping times for which each $z_{\cT_k}$ is finite and conclude with
Gronwall's inequality that $z_{\cT_k}$ is uniformly bounded over all
$k\in\N$ and $n,m\in\N$. The result follows by monotone convergence.
\end{pf}

\begin{lemma}\label{estimatesdes2}
Suppose that
\[
\sup_{t\in[0,T]} \bigl\llvert H_t^{n,m}\bigr
\rrvert %
\]
is uniformly bounded over all $n,m$ and that $Y^m=Y$ for all $m\in\N$ or
\[
\sup_{n,m\in\N} \E \Bigl[\sup_{t\in[0,T]} \bigl
\llvert \cZ _t^{n,m}\bigr\rrvert ^2\Bigr]<
\infty. %
\]
If additionally
%
\begin{equation}
\label{eq1802-2} \lim_{m\to\infty}\limsup_{n\to\infty}\E
\Bigl[\sup_{t\in
[0,T]}\bigl\llvert \cZ^{n,m}_t-
\overline\cZ^{n,m}_t\bigr\rrvert ^2 \Bigr]=0,
\end{equation}
then,
\[
\lim_{m\to\infty}\limsup_{n\to\infty}\E \Bigl[\sup
_{t\in
[0,T]}\bigl\llvert \cU^{n,m}_t-\overline\cU^{n,m}_t\bigr\rrvert ^2 \Bigr]\to0\qquad
\mbox{as }n\to\infty. %
\]
\end{lemma}

\begin{pf}
We rewrite, for $t\in[0,T]$,
\begin{eqnarray*}
\cU^{n,m}_t-\overline\cU^{n,m}_t&=&\int
_0^t \bigl(\cU^{n,m}_{s-}-
\overline\cU ^{n,m}_{s-}\bigr)H^{n,m}_s \,\dd
Y_s
-\int_0^t\cU^{n,m}_{s-}H^{n,m}_s
\,\dd \bigl(Y-Y^m\bigr)_s
\\
&&{}+\cZ^{n,m}_t
-
\overline\cZ^{n,m}_t.
\end{eqnarray*}
We fix $n,m\in\N$ and consider $z(t)=\E[\sup_{s\in[0,t]}\llvert  \cU
^{n,m}_s-\overline\cU^{n,m}_s\rrvert  ^2]$ for $t\in[0,T]$. Further, denote by
$\kappa_1$ a uniform bound for $\sup_{n,m\in\N} \sup\llvert  H_t\rrvert  ^{n,m}$
and, if applicable, for $\sup_{n,m} \E[\sup_{t\in[0,T]} \llvert   \cZ
_t^{n,m}\rrvert  ^2]$. Using that $(a_1+a_2+a_3)^2\leq3(a_1^2+a_2^2+a_3^2)$
$(a_1,a_2,a_3\in\R)$ and Lemma~\ref{estimateintegral}, we get that
\begin{eqnarray*}
z(t)&\leq&3\kappa_2\kappa_1^2\int
_0^t z(s) \,\dd s+ 3 \E \biggl[\sup
_{s\in[0,t]}\biggl\llvert \int_0^s
\cU^{n,m}_{s-}H^{n,m}_u \,\dd
\bigl(Y-Y^m\bigr)_u\biggr\rrvert ^2 \biggr]
\\
&&{} +3 \E \Bigl[\sup_{s\in[0,t]}\bigl\llvert \cZ^{n,m}_s-
\overline\cZ ^{n,m}_s\bigr\rrvert ^2 \Bigr]
\end{eqnarray*}
with $\kappa_2$ being uniformly bounded. In view of~(\ref{eq1802-2}),
the statement follows with Gronwall's inequality, once we showed that
\[
\lim_{m\to\infty}\limsup_{n\to\infty}\E \biggl[\sup
_{t\in
[0,T]}\biggl\llvert \int_0^t
\cU^{n,m}_{s-}H^{n,m}_s \,\dd
\bigl(Y-Y^m\bigr)_s\biggr\rrvert ^2
\biggr]=0. %
\]
If $Y=Y^m$, this is trivially true. In the remaining case, we can apply
Lemma \ref{estimatesde2} due to the uniform boundedness of $\E[\sup_{t\in[0,T]} \llvert  \cZ_t^{n,m}\rrvert  ^2]$ and conclude with Doob's
$L^2$-inequality and the martingale property of $Y-Y^m$ that
\begin{eqnarray*}
\E \biggl[\sup_{t\in[0,T]}\biggl\llvert \int_0^t
\cU^{n,m}_{s-}H^{n,m}_s \,\dd
\bigl(Y-Y^m\bigr)_s\biggr\rrvert ^2 \biggr]&\leq&4\int_0^T \E \bigl[\bigl\llvert \cU
^{n,m}_{s-}\bigr\rrvert ^2 \bigl\llvert
H^{n,m}_s\bigr\rrvert ^2 \bigr] \,\dd\bigl
\langle Y-Y^m\bigr\rangle_s
\\
&\leq&4 \kappa_1^2\kappa_3 T \int
_{V_m^c} x^2 \nu(\dd x)
\end{eqnarray*}
with $\kappa_3$ denoting the constant appearing in Lemma~\ref
{estimatesde2}. All constants do not depend on $n,m$ and the latter
integral tends to $0$ as $m\to\infty$.
\end{pf}

We denote by $\tau_1,\tau_2,\ldots$ adapted c\`adl\`ag processes with
$\tau_n(t)\leq t$ for all $t\in[0,T]$ and focus on perturbation
estimates for the processes $\cU^n=(\cU^n_t)_{t\in[0,T]}$ and $\overline
\cU^n=(\overline\cU^n_t)_{t\in[0,T]}$ given as solutions to
\[
\cU^n_t=\int_0^t
\cU^n_{\tau_n(s-)}H^n_{s} \,\dd
Y_s+\cZ^n_t
\]
and
\[
\overline\cU^n_t=\int_0^t
\overline\cU^n_{\tau_n(s-)}H_{s} \,\dd Y_s+\overline
\cZ^n_t. %
\]

\begin{lemma}\label{estimationsde}
1. (Stochastic convergence) If:
\begin{longlist}[(a)]
\item[(a)]$\tau_n(t)=t$, for $t\in[0,T]$,

\item[(b)] $\cZ^n-\overline\cZ^n\to0$ and $H^n- H\to0$ in ucp, as $n\to
\infty$, and
\item[(c)] the sequences $(\sup_{t\in[0,T]}\llvert  \cZ^n_t\rrvert\dvtx  n\in\N)$ and
$(\sup_{t\in[0,T]}\llvert  H^n_t\rrvert\dvtx  n\in\N)$ are tight,
\end{longlist}
then
\[
\cU^n-\overline\cU^n\to0\qquad\mbox{in ucp, as } n\to
\infty. %
\]

2. (Moment estimates) Let $p\geq2$. If:
\begin{longlist}[(a)]
\item[(a)] $Y$ has L\'evy measure $\nu$ satisfying $\int\llvert  x\rrvert  ^{p} \nu
(\dd x)<\infty$, and

\item[(b)] the expressions
\[
\sup_{t\in[0,T]}\bigl\llvert H^n_t\bigr
\rrvert \quad\mbox{and}\quad\E \Bigl[\sup_{t\in
[0,T]}\bigl\llvert
\cZ^n_t\bigr\rrvert ^{p} \Bigr]
\]
are uniformly bounded over $n\in\N$,
\end{longlist}
then
\[
\sup_{n\in\N} \E \Bigl[\sup_{t\in[0,T]}\bigl\llvert
\cU^n_t\bigr\rrvert ^{p} \Bigr]<\infty.
\]
\end{lemma}

\begin{pf}
(1) Statement~1 follows when combining Theorems 2.5(b) and 2.3(d) in~\cite{jacodprotter}.

(2) Since $\int\llvert  x\rrvert  ^{p} \nu(\dd x)<\infty$ the process
$(Y_t)$ has bounded $p$th moment and the statement can be proved
similarly as Lemma~\ref{estimatesde2} by using Lemma~\ref
{estimateintegral} and Gronwall's inequality.
\end{pf}
\end{appendix}

\section*{Acknowledgement}
 We thank two anonymous referees for their valuable comments.


%

\printaddresses
\end{document}